\documentclass[11pt]{article}

\usepackage{
		float,
		graphicx,		
		amsmath,		
		amssymb,		
		amsthm,			
		enumerate,		
		thmtools,		
		xspace,			
      		hyperref,		
      	        mathtools
              }
        \usepackage[margin=2.8cm]{geometry}	

\declaretheorem[style=definition,numberwithin=section]{Definition} 
\declaretheorem[style=plain,sibling=Definition]{Theorem}
\declaretheorem[style=plain,sibling=Definition]{Lemma}
\declaretheorem[style=plain,sibling=Definition]{Proposition}
\declaretheorem[style=plain,sibling=Definition]{Corollary}

\newcommand{\proper}[0]{proper }
\newcommand{\feas}[0]{feasible }

\newcommand{\p}[1]{\mathbb{P}\left[#1\right]} 

\newcommand{\din}[1]{{d^-_{ #1}}}
\newcommand{\dout}[1]{{d^+_{#1}}}
	\newcommand{\N}[2]{N_{{#1},{#2}}\left(n\right)}
	\newcommand{\Np}[2]{N_{{#1},{#2}}^{\per}\left(n\right)}
	\newcommand{\Nr}[2]{N_{{#1},{#2}}^{\per,r}\left(n\right)}

\newcommand{\diarray}[0]{(\diseq)_{n\in \mathbb N}}

\newcommand{\CM}[1]{\text{CM}_{n}\left( #1 \right)}

\newcommand{\instub}[1][]{{W^{-}_{#1}}}
\newcommand{\outstub}[1][]{{W^{+}_{#1}}}

\newcommand{\M}{\mathcal{M}}

\newcommand{\V}[0]{V}

\newcommand{\E}[0]{E}

\newcommand{\graphn}[0]{G_n}
 
\newcommand{\graphseq}[0]{G_{\diseq}}

\newcommand{\mgraphseq}[0]{\tilde{G}_{\diseq}}

\newcommand{\digraphspace}[0]{\mathcal{G}_{\diseq}}

\newcommand{\grapharray}{\left(\graphseq\right)_{n\in\mathbb N}}

\newcommand{\pbond}[2]{{p^{\text{bond}}_{{#1}, {#2}}}}
\newcommand{\psite}[2]{{p^{\text{site}}_{{#1}, {#2}}}}

\newcommand{\per}[1][]{\pi_{{#1}}}
\newcommand{\pcb}{{\pi_c^{\text{bond}}}}
\newcommand{\pcs}{{\pi_c^{\text{site}}}}

\newcommand{\diprobdisbond}[0]{\pbond{j}{k}}
\newcommand{\diprobdissite}[0]{\psite{j}{k}}

\newcommand{\csb}{c^{\text{bond}}}
\newcommand{\css}{c^{\text{site}}}
\newcommand{\Mp}{\mathcal{M}^{\per}}
\newcommand{\instubp}[1][]{{W^{-, \per}_{#1}}}
\newcommand{\outstubp}[1][]{{W^{+, \per}}_{#1}}

\newcommand{\instubpp}[1][]{{W^{'-, \per}_{#1}}}
\newcommand{\outstubpp}[1][]{{W^{'+, \per}}_{#1}}

\newcommand{\instubr}[1][]{{W^{-, r}_{#1}}}
\newcommand{\outstubr}[1][]{{W^{+, r}}_{#1}}
\newcommand{\instubm}[1][]{{W^{-, m}_{#1}}}
\newcommand{\outstubm}[1][]{{W^{+, m}}_{#1}}
	
\newcommand{\diseq}[1][]{{{\bf d}_{{#1}}^n}}

\newcommand{\diseqrv}[1][]{{{\bf D}_{\pi}^n}}

\newcommand{\diseqp}{\diseq[\per]}

\newcommand{\pgraph}{G_{n}^{\per}} 
\newcommand{\mpgraph}{\tilde{G}_{\diseq}^{\per}} 

\newcommand{\pgrapharray}{\left(\pgraph\right)_{n \in \mathbb N}}
\newcommand{\mpgrapharray}{\left(\mpgraph\right)_{n \in \mathbb N}}
\newcommand{\pdiarray}[0]{\left(\diseqp\right)_{n \in \mathbb N}}

\newcommand{\SCC}[1]{\text{ \rm SCC}\left( #1 \right)}

\newcommand{\Scc}[1][\graphseq]{\mathcal{C}\left({#1}\right)}

\newcommand{\diprob}[2]{p_{#1,#2}}

\newcommand{\dg}[2]{U\left( #1, #2 \right)}
\newcommand{\dgd}[2]{\dg{#1}{#2} := \sum_{j,k=0}^{\infty} \diprob{j}{k} {#1}^j {#2}^k }

\newcommand{\mum}[1][]{\mu_{{#1}}}

\newcommand{\D}[1][]{\mathcal D_{#1}}
\newcommand{\measure}{\nu}
\newcommand{\Qn}{Q_n'}
\newcommand{\Q}{Q}
\newcommand{\Xkn}{X_{\kappa,n}}
\newcommand{\Xk}{\tilde{X}_{\kappa}}

\newcommand{\si}{s^-}
\newcommand{\so}{s^+}
\newcommand{\ri}{r^-}
\newcommand{\ro}{r^+}

\newcommand{\pgraphn}{G_{n}^{\per}}

\usepackage{authblk}

\title{Percolation of a strongly connected component in simple directed random graphs with a given degree distribution}
\author[1]{Femke van Ieperen}
\author[1,*]{Ivan Kryven}
\affil[1]{\small Mathematical Institute, Utrecht University, PO Box 80010, 3508 TA Utrecht, 
the Netherlands.}
\affil[*]{\small i.kryven@uu.nl} 
\date{}
\begin{document}
\maketitle

\begin{abstract}

We study site and bond percolation on directed simple random graphs with a given degree distribution and derive the expressions for the critical value of the percolation probability above which the giant strongly connected component emerges and the fraction of vertices in this component.\\

\end{abstract}

\section{Introduction}
\label{sc:intro}
Percolation on infinite graphs is typically studied in the setting where edges (or vertices) are removed uniformly at random and  some connectivity-related property 
is being traced as a function of the percolation probability $\pi\in(0,1)$ that a randomly chosen edge (or vertex) is  present. 
Conventionally, such a traced property is chosen to describe connected components or clusters, that is maximal vertex sets in which any pair of vertices is connected with a path. 
Many results about sizes of connected components are known for percolation in infinite \emph{lattices} \cite{bollobas2006} and
 \emph{random graphs} \cite{fountoulakis2007percolation,durrett2007random,janson2009percolation,bollobas2010percolation,amini2010bootstrap}. 
Closely related to percolation are studies of random graph models that continuously depend on a parameter, such as the well-studied Erd\H{o}s R\'enyi random graph $G(n,p)$,  but also various models for directed random graphs, often refered to as $D(n,p)$ or $\vec G(n,p)$ \cite{luczak1992giant,luczak2009critical,bloznelis2012birth,goldschmidt2019scaling,pittel2016asymptotic}.

 Fountoulakis~\cite{fountoulakis2007percolation} and Jason~\cite{janson2009percolation}  studied percolation in random graphs with well-behaved degree sequences using techniques that rely on Molloy and Reed's theorem~\cite{molloy1995critical,molloy1998size}, which indicates whether a simple undirected random graph with a given degree sequence contains a giant component and how large it is.  These authors  showed that if one starts with a simple random graph generated by the configuration model and then removes edges (vertices) uniformly at random, the resulting percolated  graph can again be studied with the configuration model, albeit with a modified degree sequence. In this paper we construct a similar argument for studying percolation in directed graphs, also known as digraphs. Although there are several points of analogy, directed graphs generally require a distinct treatment from that of undirected graphs.
 
In \emph{digraphs}, there exist several non-equivalent definitions for a connected component, all of which give rise to an interesting percolation problem. 
Let $G=(V,E)$ be a simple digraph and $n=|V|$.  We say that $\mathcal C  \subset V$ is a strongly connected component (SCC) if for all $v_1,v_2\in \mathcal C  $ there are directed paths that connect $v_1$ with $v_2$ and $v_2$ with $v_1,$ and no other vertex from $V$ can be added to $\mathcal C$ without losing this property.  Suppose $G_n$ is uniformly sampled from the set of all digraphs with a fixed graphic degree sequence
\begin{align}
\label{eq:di_degseq1}
 \diseq :=
\left(
(\din{1}, \dout{1}),
(\din{2}, \dout{2}),\dots,
(\din{n}, \dout{n})
\right),
\end{align}
where $\din{v}$ and $\dout{v}$ indicate correspondingly the in- and out-degree of vertex $v\in V$.
Let additionally
$\mu:= \lim_{n\to \infty} \mu(n)$ where $\mu(n):=n^{-1}\sum_{v\in V}{ \din{v}} =n^{-1}\sum_{v\in V}{ \dout{v}},$
be the expected number of edges per vertex and $\mu_{11}:=\lim_{n\to \infty} \mu_{11}(n),$ with $ \mu_{11}(n):=n^{-1}\sum_{v\in V}{ \din{v}\dout{v}}.$
Several authors have formulated the existence criteria for the giant component in the context of directed graphs, see for example Penrose \cite{penrose2016}, Coulson \cite{coulson2019} and Cooper and Frieze \cite{cooper2004size}.
Cooper and Frieze considered sequences  $(\diseq)_{n\in \mathbb N}$ and showed that under certain regularity conditions on the degree sequences and providing
 that
 $\mum[11] - \mum>0,$
  the size of the largest strongly connected component $\Scc[G_n]$ is of the order $n$,  
 \begin{align}
 \label{eq:c}
\lim_{n\to\infty} \frac{  |\Scc[G_n]|}{n}=c>0.
\end{align}
Moreover, a SCC with this property is unique in the sense that the size of the second largest SCC is $o(n).$ 
If the latter limit holds, we say the random graph contains a giant strongly connected component (GSCC).
 Likewise, if the sign of the inequality is flipped, $\mum[11] - \mum<0,$ then the size of all SCCs is $o(n),$ and the random graph is said to contain no GSCC.
This result shows that there are two classes of limiting degree sequences, those that correspond to the size of the largest SCC being $\Theta(n)$ 
and those for which the size is $o(n)$. 

 In this paper we study the percolated graph $G_n^\pi$, in which each edge in $G_n$ is randomly removed with probability $1-\pi$. We show that if the GSCC exists in the original graph $G_n$, removing a positive fraction of edges (or vertices) can modify the degree distribution just in the right way to flip the sign of the inequality, while keeping
 the percolated graph to be uniform in the a set of all graphs with a fixed (modified) degree distribution.  
  By combining the latter observation with the theory of Cooper and Frieze, 
we show that the `phase transition' from $\Theta(n)$ to $o(n)$ takes place at a  critical value  
$\pi_c= \frac{\mum}{\mum[11]}$, such that only for $\pi>\pi_c$,  $G_n^\pi$ contains a GSCC with high probability (w.h.p.).
 The critical threshold $\pi_c$ is the same for bond and site percolation, and the expressions for $c(\pi),$ as used in \eqref{eq:c}, are closely related, $\css(\pi) = \per \csb(\pi)$.
 This work and the related proofs, are inspired by the results for percolation in undirected graphs by Fountoulakis \cite{fountoulakis2007percolation} and are based on Chapter 4 of Thesis \cite{Ieperen2020}.

\section{Main result}
\label{sc:main}
 This section introduces our main theorems for the percolation threshold of the GSCC.  We consider two types of percolation processes on simple digraph $\graphn=(V_n,E_n),$ $n=|V_n|$ that result in a random subgraph $\pgraphn$ on the same vertex set:
\begin{itemize}
	\item \textit{Bond percolation,} fix percolation probability $\per \in (0,1)$, then each edge of $\graphn$ is removed independently of the other edges with probability $1-\per$. 
	\item \textit{Site percolation,} fix percolation probability $\per \in (0,1)$, for each vertex of $\graphn$, all the edges incident to this vertex are removed with probability $1-\per$ independently of the other vertices. Such a vertex is then referred to as a \emph{deleted vertex}.
\end{itemize}
It should be clear from the context which type of percolation is discussed.  Strictly speaking, existence of the GSCC is a limiting property of a sequence of graphs $(G_n)_{n\in \mathbb N}$, in which each element is defined by a finite graphic degree sequence $\diseq$. Thus we refer to an infinite sequence of degree sequences, $(\diseq)_{n\in \mathbb N},$ as the \emph{degree progression}, where $n$ is the index and the number of vertices in the $n^\text{th}$ element of this progression.
 Although our ultimate goal is to make statements about random graphs satisfying a specific degree distribution in the limit $n\to \infty$, the bulk of the paper is spent on determining whether a growing degree progression $(\diseq)_{n\in \mathbb N}$ maintains or acquires some property of interest w.h.p. 

To make valid statements about the GSCC, we need to impose several requirements on the degree progression.  We are interested in simple graphs, which indicates that each degree sequence in the progression has to be \emph{graphic} as required by the equivalent of Erd\H{o}s-Gallai theorem for directed graphs \cite[Theorem $4$]{lamar2011directed}.
In Section~\ref{sc:random_graphs}, we progressively add several more technical constraints on $\diarray,$ namely Definitions~\ref{def:valid_di},~\ref{def:feasible_di}, and~\ref{def:proper_deg_progression}, which we jointly refer to as the requirements for a \emph{\proper degree progression}. The latter condition guarantees a sufficient regularity of the degree progression to allow us to reason about the limiting behaviour of the corresponding random graphs and connected components therein.

\begin{Definition}
	\label{def:pc}
	The \emph{percolation threshold} of a proper degree progression $\diarray$ is given by
	\begin{align}
\label{eq:per_c_def}
\per[c] = \sup\left\{ \per \in (0,1) \bigm \vert 
\forall  \varepsilon>0, \lim_{n \rightarrow \infty} \p{   \frac{\left\lvert \Scc[\pgraph]\right\rvert}{n} \geq \varepsilon       } = 0
  \right\},
\end{align}  
where superscripts are  used to further specify the type of percolation, \emph{i.e.} $\pcb$ or $\pcs.$
\end{Definition}
\noindent For each $n$, the probability in this definition is taken with respect to $G_n^\pi$ --  random graphs that remain after percolation on uniform simple random graphs obeying $(\diseq)_{n\in \mathbb N}$. 

The following theorems can be regarded as generalisation of \cite[Theorem 1.1]{fountoulakis2007percolation}  to digraphs.  They determine the percolation threshold for the existence of the giant strongly connected component, and, if this threshold exists, the theorems additionally identify the fraction of the vertices in this component for the bond and site percolation processes.
\begin{Theorem}
	\label{thm:percolation_di_strong}
	Let $\diarray$ be a proper degree progression and  $\mum[11] (n)- \mum(n)>0$ for all $n$.		 Then the thresholds for the emergence of the giant strongly connected component during bond and site percolation are the same and equal to:
$
	\pi_c = \frac{\mum}{\mum[11]}<1.
$
	\end{Theorem}
	
\noindent Let $\N{j}{k}$ be the number of vertices with in-/out-degree $(j,k)$ in  $G_n,$
and let
$
p_{j,k}:=\lim\limits_{n \rightarrow \infty}\frac{\N{j}{k}}{n}
$
exist.
Let additionally,
$$U_{\per}^\text{bond}(x,y) := \sum\limits_{j,k\geq0} p_{j,k} (1-\per+\per x)^j (1-\per+\per y)^k,$$ 
$$U_{\per}^-(x) := \left(\per\mu\right)^{-1}\frac{\partial}{\partial y}U_{\per}^\text{bond}(x,y)|_{y=1}\quad
\text{and}
\quad U_{\per}^+(y) := \left(\per\mu\right)^{-1}\frac{\partial}{\partial x}U_{\per}^\text{bond}(x,y)|_{x=1},$$ 
be formal power series in $x$ and $y$, having $\per\in(0,1)$ as a parameter.
	
\begin{Theorem}	\label{thm:percolation_di_strong2}
	   Let $\diarray$ be a proper degree progression and $\per \in( \pi_c,1),$ then there are unique values $\csb(\pi)$ and $\css(\pi)$, 
	   such that for all $\varepsilon_1, \varepsilon_2>0$, 
	   $ \lim\limits_{n \rightarrow \infty} \p{ \left|\frac{\left\lvert \Scc[\pgraph]\right\rvert}{n} -\csb (\pi)\right|\geq\varepsilon_1} = 0, $
	    for the bond percolation, and
	   $ \lim\limits_{n \rightarrow \infty} \p{\left|\frac{\left\lvert \Scc[\pgraph]\right\rvert}{n} - \css(\pi)\right|\geq \varepsilon_2} = 0, $
	   for the site percolation,
	  where
	  $$
  \csb(\pi) = 1-U_{\per}^\text{bond}(x^*,1)-U_{\per}^\text{bond}(1,y^*)+U_{\per}^\text{bond}(x^*,y^*), \;  \css(\pi) = \per \csb(\pi)
  $$
  and $x^*,y^*\in(0,1)$ are the unique solutions of
$x^*= U_{\per}^-(x^*),$ 
and
$y^*= U_{\per}^+(y^*)$ correspondingly.
\end{Theorem}

The remainder of the paper is structured as follows. Section~\ref{sc:digraphs} lays out the technical premise:
 Subsection~\ref{sc:random_graphs} introduces different classes of degree sequences. Subsection~\ref{subsc:CM}
 introduces the link between random digraphs and the directed configuration model by repurposing several theorems available for undirected graphs. Subsection~\ref{sc:giant_di} gives the definition of a giant strongly connected component and formulates the Cooper and Frieze's existence theorem, Theorem~\ref{thm:giant_scc}.
Subsection~\ref{subsc:def_theorems} formulates a corollary of McDiarmid's inequality (Corollary~\ref{cor:a-b_sets}), which is used in the proof of the main result.
Section~\ref{sc:theorem_per} proves Theorems~\ref{thm:percolation_di_strong} and~\ref{thm:percolation_di_strong2} separately for  the cases of bond and site percolation  in respectively Subsections~\ref{subsc:bond} and~\ref{subsc:site}. Both of these subsections have a similar structure: first, we show that the configuration model introduced in Section~\ref{subsc:CM} can be used to study percolated digraphs; second, we find the degree distribution after percolation; third, we show  that the corresponding degree progression is almost surely \emph{proper} in large graphs and therefore the existence  theory from Section~\ref{sc:giant_di} is applicable.

\section{Random digraphs}\label{sc:digraphs}
\subsection{Degree sequence, degree progression, and degree distribution}
\label{sc:random_graphs}
A degree sequence, as introduced in equation \eqref{eq:di_degseq1}, can be uniquely defined by adopting the lexicographic order, see Ref. \cite{lamar2011directed} for details.
 Let $\digraphspace$ be the set of \emph{all} directed multigraphs obeying degree sequence $\diseq$.
Since we are interested in sampling from $\digraphspace$, we want to be sure that for a given $\diseq$,  $\digraphspace \neq \emptyset$.
This is always the case for \emph{valid} degree sequences.
\begin{Definition}
	\label{def:valid_di}
	A degree sequence $\diseq$ is \emph{valid} if  
$
	m :=\sum_{i=1}^n \din{i} = \sum_{i=1}^n \dout{i},
$
	 where  $m=|E|$  is the number of edges.
\end{Definition}
  \noindent If $\digraphspace$ contains  a simple graph, $\diseq$ is called \emph{graphical}. Theorem $4$ in \cite{lamar2011directed} gives necessary and sufficient   criteria for $\diseq$ to be graphical.
Let $\mgraphseq\in\digraphspace$ be random directed multigraph uniformly chosen from $\digraphspace$, and  $\graphseq\in\digraphspace$ -- a uniformly chosen simple digraph.
To be on the safe side and make sure that the subset of $\digraphspace$  containing simple digraphs does not have vanishing measure for large $n$, we impose restrictions on the limiting behaviour of $\diseq$. Let 
\begin{align*}
	d_{\max}(n) := \max \left\{ \max\{\din{1}, \din{2}, \ldots, \din{n}\},  \max\{\dout{1}, \dout{2}, \ldots, \dout{n}\}\right\}
\end{align*}
be the largest degree in $\diseq$ for each index $n\in \mathbb{N}$. We will refer to this quantity as simply  $d_{\max}$, implicitly assuming dependence on $n$ throughout the paper.
\begin{Definition}
	\label{def:feasible_di}
	A  degree progression $\diarray$ is called \emph{\feas}if 1) all $\diseq$ are graphical, 2) the bivariate probability distribution 	
	 $ \diprob{j}{k}:= \lim\limits_{n \rightarrow \infty}\frac{\N{j}{k}}{n}$ has finite partial moments:
	$$
		 \mu_{il}:=\lim\limits_{n\rightarrow \infty} \sum_{j,k=0}^\infty \frac{j^ik^l\N{j}{k}}{n} =  \sum_{j,k=0}^\infty j^ik^l\diprob{j}{k}\in (0,\infty), \; i,l=0,1,2,
	$$
	and 3) $d_{\max} = \mathcal{O}\left(\sqrt{n}\right)$. Note that according to Definition~\ref{def:valid_di}, we have $\mu:=\mu_{10}=\mu_{01}$.
\end{Definition}
\noindent  We refer to $p_{j,k}$ as the \emph{degree distribution} of a feasible degree progression. 

Finally, to apply the theory developed by Cooper and Frieze for existence of the giant strongly connected component \cite{cooper2004size}, we further narrow down the class of degree progressions.\\
\begin{Definition}
	\label{def:proper_deg_progression}
	A degree progression $\diarray$ is \emph{\proper} if it is \feas and additionally satisfies
	\begin{enumerate}
		\item $d_{\max} \leq \frac{n^{1/12}}{\ln n}$,
		\item $ \rho(n): = \max \left(\sum_{j,k=0}^\infty \frac{j^2k\N{j}{k}}{\mum n}, \sum_{j,k=0}^\infty \frac{jk^2\N{j}{k}}{\mum n}\right) = o(d_{\max})$.
	\end{enumerate}	
\end{Definition}
\noindent  Thus far, we have defined a chain of classes for $\diarray$: valid $\supset$ graphical $\supset$ feasible $\supset$ proper, where the definitions of valid and graphical are extended from $\diseq$ to $\diarray$ element-wisely.

\subsection{Directed configuration model}
\label{subsc:CM} 
The behaviour of \emph{simple} random digraphs can be studied with the \emph{directed configuration model}, as defined bellow.
\begin{Definition}
	\label{def:configuration}
	Let $\diseq$ be a valid degree sequence. For all  vertices enumerated with $ i \in [n]$, let the set of \emph{in-stubs} $\instub[i]$ consist of $\din{i}$ unique elements 
	and the set \emph{out-stubs} $\outstub[i]$ contain $\dout{i}$ elements. Let $\instub = \cup_{i \in [n]} \instub[i]$ and $\outstub = \cup_{ i \in [n]}\outstub[i]$. Then a \emph{configuration}  $\M$ is a random perfect bipartite matching of $\instub$ and $\outstub$, 
that is a set of tuples $(a,b)$ such that each tuple contains one element from $\instub$ and one from $\outstub$ and each element of $\instub$ and $\outstub$ appears in exactly one tuple.
\end{Definition} 
 \noindent A configuration $\M$  prescribes a matching for all stubs, and therefore,  defines a multigraph $\mgraphseq$ with vertex set  $V = [n]$ 
 and edge multiset 
\begin{align}
\label{eq:CM_edge}
\E = [ (i,j) \mid   \outstub[i] \ni a,   \instub[j] \ni b, \text{ and }(a,b) \in \M  ].
\end{align}
 Note that multiple configurations may correspond to the same graph. We will now study the probability that the configuration model  generates a specific multigraph $\mgraphseq$. 
\begin{Proposition}
	\label{prop:CM_prob}
	Let $\mgraphseq$ be a multigraph with degree sequence $\diseq$, and   $\Upsilon_{i,j}$  be the multiplicity of edge $(i,j)$ in $\mgraphseq$. Then there holds
	\begin{align}
	\label{eq:CM_prob}
	\p{\CM{\diseq} = \mgraphseq} = \frac{1}{m!}\frac{\prod_{i=1}^n \din{i}!\prod_{i=1}^n\dout{i}!}{\prod_{1\leq i,j\leq n} \Upsilon_{i,j}!}.
	\end{align}
	\begin{proof}
		This proof is based on~\cite[Proposition $7.4$]{van2016random} formulated for undirected graphs.
 There are $m!$ different configurations. As the configuration is chosen uniformly at random, 
		\begin{align*}
		\p{\CM{\diseq} = \mgraphseq} = \frac{1}{m!}N\left(\mgraphseq\right),
		\end{align*} 
		with $N\left(\mgraphseq\right)$ being the number of distinct configurations inducing  $\mgraphseq$.  It follows from equation \eqref{eq:CM_edge} that permuting the stub labels results in a different configuration that induces the same multigraph. There are  $\prod_{i=1}^n \din{i}!\prod_{i=1}^n\dout{i}!$ such permutations. 
				 For $a,a' \in \outstub[i]$ and $ b,b' \in \instub[j]$ with $(a,b), (a',b') \in \M$, any permutation swapping $a$ with $a'$ and $b$ with $b'$ results in the same configuration. We  compensate for this by a factor $\Upsilon_{i,j}!$ to obtain:
		$$
		N\left(\mgraphseq\right) = \frac{\prod_{i=1}^n \din{i}!\prod_{i=1}^n\dout{i}!}{\prod_{1\leq i,j\leq n} \Upsilon_{ij}!},
		$$
		which completes the proof. 
	\end{proof}	
	\end{Proposition}
	
	\begin{Corollary}
	\label{prop:CM_prob_A}
	 Conditional on the event that configuration model generates a simple digraph, an element of $\digraphspace$ is chosen uniformly. 
	 \begin{proof}
 	Fix $\diseq$ and let the edge multiplicity $\Upsilon_{i,j} \leq1$ for all edges, then the sampling probability \eqref{eq:CM_prob} is constant. Hence, the model  generates  all elements $\graphseq\in \digraphspace$ with equal probability. 
	 \end{proof}
	\end{Corollary}
 Generally speaking, we are interested in the statements of the form
  \begin{equation}\label{eq:A}
  \lim_{n \rightarrow \infty} \p{ \mgraphseq \in \mathcal{A}\left(\diseq\right) } =1,
 \end{equation}
where $\mathcal{A}\left(\diseq\right)$ is a set of all graphs that obey the degree sequence $\diseq$ and additionally satisfy some desired property. 
  If this limit holds for a given property $\mathcal{A}$, then we say that the random graph has this property with high probability (w.h.p) or asymptotically almost surely (a.a.s.). The goal of the remainder of this section is to show that if equation \eqref{eq:A} holds, then 
\begin{align*}
\lim_{n \rightarrow \infty}  \p{\mgraphseq \in \mathcal{A}\left(\diseq\right) \mid \mgraphseq \,\text{is simple}} =1.
\end{align*} 
First,  we show that the probability that the configuration model generates a simple graph is bounded away from zero. 
	\begin{Proposition}\cite[Theorem $4.3$]{chen2013directed}
		 \label{thm:di_simple}
		Let $\diarray$ be a \feas degree progression. 
		 The probability that the configuration model generates a simple graph is asymptotically 
	$
	 e^{-\frac{\mum[11]}{\mum}  - \frac{\left(\mum[20] -\mum\right)\left(\mum[02]-\mum\right)}{\mum} } >0.
	$
	\begin{proof}
		The proof follows from the proof of~\cite[Theorem $4.3$]{chen2013directed}. 
		 The main difference is that in~\cite{chen2013directed} the in-degree of a vertex is independent of its out-degree. It suffices to replace Condition $4.1$ and Lemma $5.2$ from Ref.~\cite{chen2013directed} with the requirement of a \feas degree progression. 
	\end{proof}
\end{Proposition}
	\begin{Lemma}
	\label{lm:multi_simple_0}
	Let $\diarray$ be a \feas degree progression, 
	and let $\mathcal{A}\left(\diseq \right)$ be a set of multigraphs all satisfying  $\diseq$.   Let  $\mgraphseq$ be a random multigraph generated by the configuration model.\\1)
	 If
		$
	\lim\limits_{n \rightarrow \infty}  \p{ \mgraphseq \in \mathcal{A}\left(\diseq\right)} =0,
	$ then
	$
	\lim\limits_{n \rightarrow \infty}  \p{\mgraphseq \in \mathcal{A}\left(\diseq\right) \mid \mgraphseq \,\textrm{\rm is simple}} =0.
	$\\
	2) If $
	\lim\limits_{n \rightarrow \infty}  \p{ \mgraphseq \in \mathcal{A}\left(\diseq\right)} =1,
	$ then 
	$
	\lim\limits_{n \rightarrow \infty}  \p{\mgraphseq \in \mathcal{A}\left(\diseq\right) \big|\, \mgraphseq \,\textrm{\rm is simple}} =1.
	$
	\begin{proof}1)
		By Bayes' rule 
		$
		\p{\mgraphseq \in \mathcal{A}\left(\diseq\right) \mid \mgraphseq \,\text{is simple}}  \leq \frac{\p{\mgraphseq \in \mathcal{A}\left(\diseq\right)}}{\p{ \mgraphseq\, \text{is simple}}}.
		$
		Proposition~\ref{thm:di_simple} assures that the denominator does not vanish: 
$\p{ \mgraphseq\, \text{is simple}} = \left(1+o(1)\right) e^{-\frac{\mum[11]}{\mum}  - \frac{\left(\mum[20] -\mum\right)\left(\mum[02]-\mum\right)}{\mum} },$
hence
		$
		\liminf\limits_{n\rightarrow \infty}  \p{ \mgraphseq\, \text{is simple}} >0;
		$
		while the numerator converges to zero by assumption. \\
		2) It is enough to apply 1) to the complement event $\overline{\mathcal{A}}\left(\diseq\right) := \digraphspace\setminus \mathcal{A}\left(\diseq\right)$.		
	\end{proof}
\end{Lemma}

\subsection{Giant strongly connected component in a directed graph}
\label{sc:giant_di}
 A natural definition of a path in a directed graph requires that a path respects edge directions:
 \begin{Definition}
 \label{def:paths}
Let $G=(V,E)$ be a digraph.  A pair of vertices $v_1, v_k \in V$  is connected by a \emph{directed path} if there exist distinct vertices $v_2,v_3, \ldots ,v_{k-1} \in \V$ such that  for all $ i \in\{2,3,\ldots, k\}$ $\left(v_{i-1}, v_i\right) \in \E$. We refer to such a sequence as a directed $v_1-v_k$ path.  
\end{Definition}
 \noindent This definition of connectivity can be extended to define connected components:
 \begin{Definition}(Strongly connected component) \label{def:scc}
	Consider a directed graph  $G$. The \emph{strongly connected components} of $G$ are the maximal subsets of $\V$ such that between any pair of vertices $u,v$, directed $u-v$ and $v-u$ paths exist simultaneously. 
\end{Definition}
\begin{Definition}\label{def:in_out_CC} 	Consider a directed graph  $G=(V,E)$ and take $v \in \V$. Then the \emph{strong-component} of $v$, denoted by $\SCC{v}$, consists of $v$ itself and all vertices $w \in \V$ for which  both  $v-w$ and $w-v$ paths exist. 
\end{Definition}
\noindent 
\noindent Let $\Scc[\graphn]$ be the largest strongly connected component.
The notion of a \emph{giant component} is introduced as a limiting property of the sequence of $\Scc[\graphn]$.
\begin{Definition}
The graph progression $\grapharray$ is said to contain a \emph{giant strongly connected component} (GSCC) if:
a) $\lim\limits_{n\rightarrow \infty} \frac{  |\Scc |}{n} = \zeta >0$ or
b) $\lim\limits_{n\rightarrow \infty}\frac{|\{v \,\vert \,  \SCC{v} \subset \Scc  \}|}{n}= \zeta >0$.\\
Criteria  a) and b) are equivalent.
\end{Definition} 
Let $x^*,y^*\in(0,1)$  be the unique solutions of respectively
$
x^*= U^-(x^*)
$
and
$
y^*= U^+(y^*)
$,
where 
$$
\begin{aligned}
U^+(y) := \mu^{-1}\frac{\partial}{\partial x}U(x,y)|_{x=1},\\
U^-(x) := \mu^{-1}\frac{\partial}{\partial y}U(x,y)|_{y=1},
\end{aligned}
$$
and $\dgd{x}{y}$
is the generating function of $p_{j,k}$. Let further,
  $$
\zeta := 1-U(x^*,1)-U(1,y^*)+U(x^*,y^*).
  $$
\begin{Theorem}[Existence of
GSCC]\cite[Theorems $1$ and $2$]{cooper2004size} 
	\label{thm:giant_scc}
	Consider a \proper degree progression $\diarray$ and a uniformly random sequence of simple graphs $\grapharray$. Then,
	\begin{enumerate}
		\item If $\frac{\mum[11]}{\mum} < 1$, 
		 the size of any SCC is $\mathcal{O}\left(d_{\max}^2\ln n\right)$ with high probability.
		\item If  $U^-(0)>0,\; U^+(0)> 0$ and $\frac{\mum[11]}{\mum} >1$ then with high probability there is a  unique giant strongly connected component with vertex set of size $\zeta n$. 
	\end{enumerate}
\end{Theorem}

\subsection{Concentration inequalities} 
\label{subsc:def_theorems}
We introduce concentration inequalities that are later used in the proof of Theorem~\ref{thm:percolation_di_strong} to treat two sources of randomness: the multigraph is random and percolation randomly removes edges.   
\begin{Theorem}
	(Hoeffding's inequality)\cite{boucheron2013concentration} \label{thm:hoeffding}
	Let $X_1,X_2, \ldots, X_n$ be independent random variables. Suppose that $a_i \leq X_i \leq b_i$ for all $i \in \{1,2,\ldots, n\}$ and define $c_i = b_i - a_i$. Furthermore define $S_n = \sum_{i=1}^n X_i$. Then there holds 
	\begin{align}
	\p{\lvert S_n - \mathbb{E}\left[S_n\right] \rvert >t} \leq 2\exp\left(-\frac{2t^2}{\sum_{i=1}^n c_i^2}\right).
	\end{align} 
\end{Theorem}
\noindent The following concentration inequality is a corollary of a theorem by McDiarmid. 
\begin{Theorem} \cite[Theorem $7.4$]{mcdiarmid1989method}
	\label{thm:mcdiarmid}
	Let $\left(V,d\right)$ be a finite metric space. Suppose there exists a sequence $\mathcal{P}_0, \mathcal{P}_1, \ldots, \mathcal{P}_s$ of increasingly refined partitions with $\mathcal{P}_0$ the trivial partition consisting of $V$ and $\mathcal{P}_s$ the partition where each element of $V$ is a partition element on its own. Take a sequence of positive integers $c_0, c_1, \ldots, c_s$ such that for all $k \in \{1,2,\ldots, s\}$ and any $A,B \in \mathcal{P}_k$ with $C$ satisfying $A,B\subset C \in \mathcal{P}_{k-1}$ there exists a bijection $\phi: A \rightarrow B$ with $d\left(x,\phi(x)\right) \leq c_k$ for all $x \in A$. Let the function $f: V \rightarrow \mathbb{R}$  satisfy $\lvert f(x) - f(y) \rvert \leq d(x,y)$ for all $x,y,\in V$. Then for $X$ uniformly distributed over $V$ and any $t>0$ there holds
	\begin{align*}
	\p{\lvert f(X) - \mathbb{E}\left[f(X)\right] \rvert > t } \leq 2 \exp\left(-\frac{2t^2}{\sum_{k=0}^s c_k^2}\right).
	\end{align*}
\end{Theorem}

\begin{Corollary}
	\label{cor:a-b_sets}
	Consider two finite sets $A_0$ and $A_1$ with $\lvert A_0 \rvert = a_0$ and $\lvert A_1\rvert = a_1$. Let $ S:= \cup_{i \in \{0,1\}} \{(x,i) \mid x \in A_i\}$. A subset of $S$ containing $b_0$ elements with $i=0$ and $b_1$ elements with $i=1$ is called a $\left(b_0,b_1\right)$-subset of $S$. Let $V$ be the space of all $\left(b_0,b_1\right)$-subsets of $S$. Let $f: V \rightarrow \mathbb{R}$ be a function such that for any $B,B' \in V$ there holds $\lvert f\left(B\right) - f\left(B'\right)\rvert \leq \lvert B \triangle B' \rvert $. Here $B \triangle B'$ denotes the symmetric difference, i.e. $B \triangle B' = \left(B\cup B'\right) \setminus \left(B \cap B'\right)$.
	Then for  $X$ distributed uniformly over $V$ and any $t > 0$ there holds
	\begin{align}
	\label{eq:cor_mcdiarmid}
	\p{ \left\lvert f(X) - \mathbb{E}\left[f(X)\right] \right\rvert > t} \leq 2\exp\left(-\frac{t^2}{8(b_0+b_1)}\right).
	\end{align}
	\begin{proof}
		Consider a $\left(b_0,b_1\right)$-subset of $S$. Assign each element a unique number from the index set $\{1,2,\ldots, b_0+b_1\}$, such that for all elements with $i=0$ this number is smaller than $b_0+1$. Note that this implies that for each element with $i=1$ its index is larger than $b_0$. A $\left(b_0,b_1\right)$-subset of $S$ with such a numbering is called  a $\left(b_0,b_1\right)$-ordering of $S$. 
		Define $W$ to be the set of all $\left(b_0,b_1\right)$-orderings of $S$.
	The function $f: V\rightarrow \mathbb{R}$ can be extended to a function $f:W \rightarrow \mathbb{R}$ by regarding each $\left(b_0,b_1\right)$-ordering as $\left(b_0,b_1\right)$-subset. This extension respects the relation $\lvert f(x) - f(y)\rvert \leq x \triangle y$, \emph{i.e.} it holds for $x,y \in W$  as well. This is true since for any two orderings $x,y$ their symmetric difference as $\left(b_0,b_1\right)$-orderings is bounded from bellow by their symmetric difference as $\left(b_0,b_1\right)$-subsets. The next step in proving equation \eqref{eq:cor_mcdiarmid} is applying  Theorem~\ref{thm:mcdiarmid} to the metric space $\left(W, \triangle \right)$. \\

We will now define a sequence of refined partitions on $W$  using the notion of an $i$-prefix. An $i$-prefix determines the first $i$ elements of an ordering. This allows for all  $k \in \{0,1,\ldots, b_0+b_1\}$ to construct the partition  $\mathcal{P}_k$ by defining its elements to be  the sets of orderings with the same $k$-prefix.  The partition $\mathcal{P}_0$ is the trivial partition consisting of $W.$ As each $\left(b_0,b_1\right)$-ordering has $b_0+b_1$ elements,  $\mathcal{P}_{b_0+b_1}$ will be the partition where each element is a single ordering. Next the values $c_k$ need to be determined. Take $B,D \in \mathcal{P}_k$ with $C$ satisfying $B,D\subset C \in \mathcal{P}_{k-1}$. This implies that any ordering in $B$ has the same $k-1$-prefix as an ordering in $D$. Furthermore these orderings must differ at the $k^\text{th}$ element. The remaining $b_0+b_1-k$ elements can be any element that is not present in the $k$-prefix that lead to a valid ordering. Denote the $k^\text{th}$ element of any ordering in $B$ by $a_{B,k}$. Similarly let $a_{D,k}$ denote the $k^\text{th}$ element of any ordering in $D$. Define the bijection $\phi: B \rightarrow D$ by  taking $x \in B$ and mapping its  $k^{th}$ element to $a_{D,k}$. If $x$ contains $a_{D,k}$ at some position $l >k$, map the $l^\text{th}$ element of $x$ to $a_{B,k}$. All the other elements are unchanged by the bijection. By definition this is an element of $D$. 

		According to definition of $\phi$ for any $x \in B$ we have $\lvert x \triangle \phi(x) \rvert \leq 4$. Thus we may take $c_k = 4$ for all $k \in \{1,2,\ldots , b_0+b_1\}$. Applying Theorem~\ref{thm:mcdiarmid} we find that equation \eqref{eq:cor_mcdiarmid} holds for any $t>0$ and  $X$ distributed uniformly over $W.$
	Remark that each element in $V$ gives rise to $b_0!+b_1!$ different orderings. All these orderings  have the same value under $f$. Thus the probability that $f(X) = c$ does not change when we take $X$ to be a uniformly random element of $V$ instead of $W$. Together with the above equation, this proves the claim.
	\end{proof}
\end{Corollary}

\section{Proofs of Theorems~\ref{thm:percolation_di_strong} and~\ref{thm:percolation_di_strong2}}
\label{sc:theorem_per}
Theorems ~\ref{thm:percolation_di_strong} and~\ref{thm:percolation_di_strong2} are proven separately for bond and site percolation by using similar techniques as in the proof of \cite[Theorem 1.1]{fountoulakis2007percolation}, which determines the percolation threshold in undirected graphs. Although these theorems are formulated  for simple digraphs, we  instead prove the statements on random directed multigraphs introduced in Section~\ref{subsc:CM}. The results on the multigraphs are then transferred to  simple digraphs using Lemma~\ref{lm:multi_simple_0}. 
After percolation with parameter $\pi$, a fixed degree sequence $\diseq$ becomes a random variable with a distribution of outcomes.
We use $\diseqrv$ to denote this random variable and  $\diseqp$ to refer to a fixed outcome of $\diseqrv$.

\subsection{Bond percolation}
\label{subsc:bond}
 Bond percolation removes edges in a graph, thus in a configuration, it removes  in-stubs  together with their matched out-stubs. Let  $\instubp$ and $\outstubp$ denote the in-stubs and out-stubs surviving percolation. Conditional on $\diseqrv=\diseqp$ there is a bijection 
 between the configuration on $\left(\instubp,\outstubp\right)$ and the configuration on $\left(\instub[\diseqp], \outstub[\diseqp]\right)$, \textit{i.e.} the stubs inducing the degree sequence $\diseqp$.   Let us fix such a bijection, and let $\D[n]$ be the probability space containing all degree sequences $\diseqp$ that can be obtained by applying percolation to a random configuration on $\left(\instub, \outstub\right)$. The probability assigned to each  $\diseqp$ is the probability that it is induced by $\left(\instubp,\outstubp\right)$. The probability space for the degree progression $\pdiarray$ is the product space  $\D = \prod_{n=1}^\infty \D[n]$ with the product measure $\measure$.

 The remainder of this section is as follows:  in Section~\ref{subsubsc:random_bond} we show  that conditional on the degree sequence after percolation, each configuration on $\instub[\diseqp]$ and $\outstub[\diseqp]$ is equally likely.   In Section~\ref{subsubsc:degree_dis_bond} we determine the limit of the expected number of vertices with degree $(j,k)$ after percolation. Combining these results in Section 
\ref{subsubsc:theorem_bond}, the proofs  of Theorems~\ref{thm:percolation_di_strong}~and~\ref{thm:percolation_di_strong2}  are completed by showing that an element of $\D$ is $\measure$-a.s.~proper, which authorises applying  Theorem~\ref{thm:giant_scc}.

 \subsubsection{A percolated configuration is a uniformly random configuration}
\label{subsubsc:random_bond}
  In this section we show that conditional on the degree sequence after percolation, the configuration on $\left(\instubp,\outstubp\right)$ is also a uniformly random one. The proof is split into two lemma's.
\begin{Lemma}
	\label{lm:configuration_points_percolation} Let  $l$ out of $m$ edges survive bond percolation applied to a uniformly random configuration $\M$ on $\left(\instub, \outstub\right).$ Then the surviving stubs $\instubp\subset \instub$ and $\outstubp \subset \outstub$ are uniformly distributed amongst all pairs of subsets of $\instub$ and $\outstub$ of size $l$. 
	\begin{proof}
		Since the graph contains $m$ matches of which $l=\lvert \instubp \rvert = \lvert \outstubp \rvert $ survive percolation,  the probability that exactly those $l$ matches remain is  $\frac{1}{{m \choose l}}.$		 
		 		It is left to investigate the probability that all stubs in $\instubp$ have their match in $\outstubp$, that is that $\M$ can be decomposed into a perfect bipartite matching of $\instubp$ with $\outstubp$ and a perfect bipartite matching of  $\instub \setminus \instubp$ with $\outstub\setminus \outstubp$. Between two sets of size $l$ there are $l!$ perfect bipartite matchings, hence the probability that $\M$ decomposes as desired, is
		$l!(m-l)! / m!$. 
		Thus the probability that  $\left(\instubp, \outstubp \right)$ are the stubs surviving percolation is
$
		\frac{l!\, (m-l)!}{m!}\frac{1}{{m\choose l}} = \frac{1}{{m\choose l}^2}.
$
		This is the probability that $\instubp\!\subset\instub$ and $\outstubp\!\subset\outstub$ are uniformly random subsets both of size $l$.	\end{proof} 
\end{Lemma}

\begin{Lemma}
	\label{lm:percolation_is_CM}
	Conditional on having degree sequence $\diseqp$ after bond percolation, i.e. $\diseqrv = \diseqp$, all configurations of $\instub[\diseqp]$ with $\outstub[\diseqp]$ are equally likely.
	\begin{proof}
	 It suffices to show that for any perfect bipartite matching $\Mp$ of $\instub[\diseqp]$ with $\outstub[\diseqp]$:
		\begin{align}
		\label{eq:1/l}
		\p{\Mp | \diseqrv = \diseqp} = \frac{1}{l!},
		\end{align}
		where $l$ is the sum of the in-degrees of $\diseqp$.
		Using $\p{\lvert \instubp \rvert =l  \middle| \diseqrv = \diseqp} =1$ and Bayes' formula we obtain:
		\begin{align}
		\label{eq:m_d_l}
		\p{\Mp | \diseqrv = \diseqp} 
		&= \p{\Mp \middle|\, \lvert\instubp\rvert=l,  \diseqrv = \diseqp}=\frac{\p{\Mp \cap \diseqrv = \diseqp  \,\middle|\, \lvert \instubp \rvert=l}}{\p{\diseqrv = \diseqp \middle \lvert \instubp \rvert=l}},
		\end{align}		
		which will  show to be equal $\frac{1}{l!}$. Let
		$S(\diseqp)$  be the collection of pairs of subsets of $\left(\instub, \outstub\right)$ that induce the degree sequence $\diseqp$.	Recalling that $\lvert\instubp\rvert = \lvert\outstubp\rvert$, we see that for any pair of subsets in $S(\diseqp)$, both sets must contain $l$ elements. In combination with Lemma~\ref{lm:configuration_points_percolation} this implies that
\begin{align}
		\label{eq:Sml2}
		\p{\diseqrv = \diseqp \middle| \, \lvert \instubp \rvert = l} =  \frac{|S(\diseqp)|}{{m\choose l}^2}.
\end{align}	
		Next, we investigate $\p{\Mp \cap \diseqrv = \diseqp  \,\middle|\, \lvert \instubp \rvert=l}$.
		 By definition, $\diseqrv = \diseqp $ implies that 
		 $\left(\instubp,\outstubp\right)\in S(\diseqp)$. Let $\left(\instubp,\outstubp\right)\in S(\diseqp)$, we aim to find the probability that the configuration on these stubs induces the configuration $\Mp$ on $\left(\instub[\diseqp],\outstub[\diseqp]\right)$. As we fixed a bijection between $\left(\instubp,\outstubp\right)$ and $\left(\instub[\diseqp],\outstub[\diseqp]\right)$, exactly one configuration of $\left(\instubp,\outstubp\right)$ induces the configuration $\Mp$ on $\left(\instub[\diseqp],\outstub[\diseqp]\right)$. 
		However we are free to choose  the configuration  on $\left( \instub \setminus \instubp, \outstub \setminus \outstubp\right)$. Thus assuming that $\left(\instubp,\outstubp\right)\in S(\diseqp),$ the probability that it induces the configuration $\Mp$ on  $\left(\instubp,\outstubp\right)$ is 
$
		\frac{(m-l)!}{m!}.
$
		Thus for any collection of remaining stubs that induce the right degree sequence, it has probability  $\frac{(m-l)!}{m!}$ to induce the right matching. However, as we condition only on the size of $\instubp$, we also must take into account the probability that exactly the desired $l$ edges survive percolation. Hence, we have:
		\begin{align*}
		\p{\Mp \cap \diseqrv = \diseqp  \,\middle|\, \lvert \instubp \rvert=l} = \frac{(m-l)!}{m!}\frac{|S(\diseqp)|}{{m\choose l}}.
		\end{align*}
		which, when combined with \eqref{eq:Sml2} and plugged into equation \eqref{eq:m_d_l}, recovers equation \eqref{eq:1/l}.
	\end{proof} 
\end{Lemma}

\subsubsection{The expected number of vertices with degree $(j,k)$ after bond percolation}
\label{subsubsc:degree_dis_bond}
In this section we derive an expression for the degree distribution after bond percolation.
\begin{Lemma}
	\label{lm:Nkj_pi}
	Let $\Np{j}{k}$ be the number of vertices in the bond-percolated graph (or configuration) with in-degree $j$ and out-degree $k$. The following limit exists
	\begin{align}
\label{eq:lim_p_per}
\pbond{j}{k}:=\lim_{n\rightarrow \infty} \frac{\mathbb{E}\left[\Np{j}{k}\right]}{n}, \text{ for } j,k\in\mathbb N_0,
\end{align}
and
\begin{align}
\label{eq:pbond}
\pbond{j}{k}  = \sum_{d^- = j}^\infty \sum_{d^+ = k}^\infty \diprob{d^-}{d^+}  {d^- \choose j}{d^+ \choose k}\per^{j+k}\left(1-\per\right)^{d^- - j + d^+ -k}.
\end{align}

	\begin{proof}
 When $j,k > d_{\max}$, $\pbond{j}{k}=\Np{j}{k}=\N{j}{k} = 0$. 
Let $j,k \leq d_{\max}$ and remark that 
\begin{align}
\label{eq:exp_np_bond}
\mathbb{E}\left[\Np{j}{k}\right] = \sum_{l=0}^m \mathbb{E}\left[\Np{j}{k}\,\middle| \, \lvert \instubp \rvert = l\right]\p{\lvert \instubp \rvert = l},
\end{align}	
which, when conditioned on $\lvert\instubp \rvert = l$, can be rewritten as  
\begin{align*}
\mathbb{E}\left[\Np{j}{k}\,\middle| \, \lvert \instubp \rvert = l\right] = \sum_{d^- =j}^{d_{\max}} \sum_{d^+ = k}^{d_{\max}} \N{d^-}{d^+}\p{ \text{$(d^-,d^+) \to (j, k)$}\,\middle| \lvert \instubp \rvert = l},
\end{align*}
where $(d^-,d^+)$ is the degree of a vertex before percolation and $(j,k)$ -- the degree after percolation. 
As Lemma~\ref{lm:configuration_points_percolation} implies that the surviving stubs are chosen uniformly at random, conditional on the size of $\instubp$, there holds
\begin{align*}
\p{ \text{$(d^-,d^+) \to(j, k)$}\,\middle| \lvert \instubp \rvert = l}= {d^- \choose j}\frac{{{m-d^-} \choose {l-j}}}{{m \choose l}} {d^+ \choose k} \frac{{{m-d^+} \choose {l-k}}}{{m \choose l}}. 
\end{align*} 
Since the edges are removed independently of each other, the size of $\instubp$ is the sum of $m$ independent Bernoulli variables, each having expectation $\per$. Applying Hoefdding's inequality (Theorem~\ref{thm:hoeffding}) shows that $\lvert \instubp \rvert$ concentrates,
\begin{align}
\label{eq:l_in_I}
\p{ \lvert \instubp \rvert  \in I} \leq \exp\left[-\Omega (\ln^2n)\right],
\end{align} 
where
\begin{equation}\label{eq:I}
I := \left[ m\per - \ln n\sqrt{n}, m\per + \ln n\sqrt{n}\right].
\end{equation}
Note that as an immediate consequence of this inequality we have $\p{l \notin I} = o\left(n^\alpha\right)$ for any $\alpha<0$.
Since we consider \proper degree progressions, there holds $d_{\max} \leq \frac{n^{1/12}}{\ln n}$, and for  $l \in I$:
\begin{align*}
\mathbb{P}\left[ \text{ $(d^-,d^+) \to(j,k)$}\,\middle| \,\lvert \instubp\rvert = l\right]
= {{d^-}\choose{j}} {{d^+}\choose{k}} \per^{j+k}\left(1-\per\right)^{d^- + d^+ - j - k}\left(1 + \mathcal{O}\left(\frac{\ln n}{n^{7/18}}\right)\right),
\end{align*}
for all $d^-,d^+ \leq d_{\max}$. In combination with equation \eqref{eq:exp_np_bond}  and the bound $\Np{j}{k} \leq n$, we find:
\begin{align}
\label{eq:bound_ex_Ndd_bond}
\mathbb{E}\left[\Np{j}{k}\right] = \left(1 +\mathcal{O}\left(\frac{\ln n}{n^{7/18}}\right)\right)\sum_{d^-=j}^{d_{\max}} \sum_{d^+=k}^{d_{\max}} \N{d^-}{d^+} {{d^-}\choose{j}} {{d^+}\choose{k}} \per^{j+k}\left(1-\per\right)^{d^- + d^+ - j - k}  +o\left(\frac{1}{n^3}\right),
\end{align}
where we set $\alpha=5$ in the estimate for $\p{l \notin I} $ when calculating the error term.
We will now show that for all  $\epsilon >0$  there exist $\kappa\left(\epsilon\right)$ and $N\left(\epsilon\right)$ such that for all $n > N,$
\begin{align}
\label{eq:bound_Ndd_bond}
		\frac{1}{n} \sum_{\mathclap{\substack{(d^-,d^+) = (0,0)\\
			d^- \geq \kappa +1 \,\text{or}\, d^+ \geq \kappa +1}}}^{(d_{\max}, d_{\max})} \N{d^-}{d^+} {{d^-} \choose j} {{d^+} \choose k} \per^{j+k}\left(1-\per\right)^{d^- + d^+ -j -k} \leq \frac{1}{n} \sum_{\mathclap{\substack{(d^-,d^+) = (0,0)\\
			d^- \geq \kappa +1 \,\text{or}\, d^+ \geq \kappa +1}}}^{(d_{\max}, d_{\max})} \N{d^-}{d^+} < \epsilon
.
\end{align}
The left inequality follows from the binomial theorem, which implies that $\sum_{j=0}^{d^-}{{d^-}\choose{j}}\per^j(1-\per)^{d^- -j} = \sum_{k=0}^{d^+}{{d^+}\choose{k}}\per^k(1-\per)^{d^+ -k} = 1$, yielding  that $ {{d^-} \choose j} {{d^+} \choose k} \per^{j+k}\left(1-\per\right)^{d^- + d^+ -j -k} \leq 1$ for all $j \leq d^-, k\leq d^+$. The right inequality holds since  $\lim\limits_{n\rightarrow \infty}\frac{\N{j}{k}}{n} = \diprob{j}{k}$ for $j,k \geq 0$, which follows from the degree progression being proper.  When combined together, equations \eqref{eq:bound_ex_Ndd_bond} and \eqref{eq:bound_Ndd_bond} prove the claim.
\end{proof}
\end{Lemma}
 Elementary calculations show that $\pbond{j}{k}$ is normalised and has the following moments:
\begin{equation}
\begin{aligned}
\label{eq:mum_per_mum11_per}
\mum^{\per, \text{bond}}:=&\sum_{j=0}^\infty \sum_{k=0}^\infty j\pbond{j}{k} = \sum_{j=0}^\infty \sum_{k=0}^\infty k\pbond{j}{k} =\per \mum,\\
 \mum[11]^{\per, \text{bond}}:=&\sum_{j,k=0}^\infty jk\pbond{j}{k} = \per^2\sum_{d^-=0}^\infty  \sum_{d^+=0}^\infty d^{-}d^{+}\diprob{d^-}{d^+} = \per^2 \mum[11].
\end{aligned}
\end{equation}
The generating function $U_{\per}^\text{bond}$ of $\pbond{j}{k}$ and auxiliary functions $U^+_{\per},U^-_{\per}$ analogous to those used in Section~\ref{sc:giant_di} are given by:
\begin{equation}
\label{eq:gfs}
\begin{aligned}
 &U_{\per}^\text{bond}(x,y) :=    U(1-\per+\per x,1-\per+\per y),\\
& U^+_{\per}(x):=U^+(1-\per+\per x),\\
&U^-_{\per}(y):=U^-(1-\per+\per y).
 \end{aligned}
 \end{equation}
Note that $U^-_{\per}(0)>0,$ $U^+_{\per}(0)>0$ for all $0<\pi<1$.
This observation, the moments \eqref{eq:mum_per_mum11_per} and generating functions \eqref{eq:gfs} provide sufficient information for applying Theorem~\ref{thm:giant_scc} to $\diprobdisbond$ and hence formally defining the percolation threshold as such a value of $\per=\hat{\per}^{\text{bond}}$  that $\mum^{\per, \text{bond}}= \mum[11]^{\per, \text{bond}}$, that is
\begin{align*}
\hat{\per}^{\text{bond}} = \frac{\mum}{\mum[11]}.
\end{align*}
It is left to show that $\diprobdisbond$ is indeed $\measure$-a.s the degree distribution of the percolated graph, and therefore, $\pcb=\hat{\per}^{\text{bond}},$ which we do in the following section.

\subsubsection{Determining $\pcb$ and $\csb$}
\label{subsubsc:theorem_bond}
Theorem ~\ref{thm:giant_scc} requires   $(\diseq)_{n\in\mathbb{N}}$ to be  a \emph{\proper}degree progression, which is assumed throughout this section. We make use of Lemma~\ref{lm:percolation_is_CM}, stating that instead of actually removing  edges, one may view the percolated multigraph $\mpgraph$ as a uniformly random \emph{configuration} obeying the percolated degree sequence $\diseqp$. 
 We show that $\pdiarray$ is indeed $\measure$-a.s. feasible (Lemma~\ref{lm:bondfeasible}) and proper (Lemma~\ref{lm:bondproper}) as an element of $\D$, and  hence, Theorem~\ref{thm:giant_scc} is applicable to $\mpgrapharray$.
  This means that Theorem~\ref{thm:giant_scc} may be applied to almost all degree progressions $\pdiarray$ to determine $\pcb$ and $\csb$ for random multigraphs. 
  Conditioning on the graph before percolation being simple will ensure that the percolated graph is simple as well. 
  Finally we apply a variant of Lemma~\ref{lm:multi_simple_0}  and show that similar assertions hold for a graph progression $\pgrapharray$ for percolated \emph{simple} graphs obeying $\diarray$, hence proving Theorem~\ref{thm:percolation_di_strong2} for the case of bond percolation.   In the remainder of this section, we make the above-stated argument formal.

\begin{Lemma}\label{lm:bondfeasible}
The degree progression after bond percolation is feasible $\measure$-a.s.
\begin{proof} We need to show that $\pdiarray$ satisfies the requirements of Definition~\ref{def:feasible_di} $\measure$-a.s.
To demonstrate that
\begin{align}
\label{eq:con_dis_bond}
\lim_{n\rightarrow \infty} \frac{\Np{j}{k}}{n} = \pbond{j}{k},\; \text{$\measure$-a.s. for } j,k\in \mathbb N_0,
\end{align}
 it suffices to show (see e.g. \cite[Lemma $6.8$]{petrov1995limit}) that for all $\epsilon >0,$ 
\begin{align}
\label{eq:petrov}
	\sum_{n=1}^\infty \p{\left\lvert \frac{1}{n}\Np{j}{k} - \pbond{j}{k}\right\rvert > \epsilon} <\infty.
\end{align}	
By definition of $\pbond{j}{k},$ for any fixed $\epsilon >0$ there is  $K$ such that for all $n > K$
\begin{align*}
\left\lvert \frac{1}{n}\mathbb{E}\left[\Np{j}{k}\right] - \pbond{j}{k} \right\rvert \leq \frac{\epsilon}{2}. 
\end{align*}
This implies that
\begin{align*}
\p{\left\lvert \frac{1}{n}\Np{j}{k} - \pbond{j}{k}\right\rvert > \epsilon} \leq \p{\frac{1}{n}\left\lvert \Np{j}{k} - \mathbb{E}\left[\Np{j}{k}\right] \right \rvert > \frac{\epsilon}{2}}.
\end{align*}
 Lemma~\ref{lm:configuration_points_percolation} states that conditional on $\lvert \instubp \rvert =l$, the stubs surviving percolation $\left(\instubp, \outstubp\right)$ are uniformly distributed amongst all pairs of subsets of $\left(\instub,\outstub\right)$ of size $l$.  $\Np{j}{k}$ is a function of $\instubp \cup \outstubp$. Furthermore for two sets $\instubp \cup \outstubp$ and $\instubpp \cup \outstubpp$ their values of $\Np{j}{k}$ differ by at most the number of elements in the symmetric difference of  $\instubp \cup \outstubp$ and $\instubpp\cup \outstubpp$. This implies that the requirements of Corollary~\ref{cor:a-b_sets} are fulfilled by setting $A_0 = \instub, A_1= \outstub, b_0=b_1=l$ and $\Np{j}{k}$ as function $f$. Applying this corollary gives:
\begin{align*}
\p{\left\lvert \Np{j}{k} - \mathbb{E}\left[\Np{j}{k}\right] \right\rvert > \frac{n\epsilon}{2} \,\vert\, \lvert \instubp\rvert = l} \leq 2\exp\left(\frac{{\epsilon}^2n^2}{64l^2}\right).
\end{align*}
For $l \in I$, defined in \eqref{eq:I}, this probability is $o\left(\frac{1}{n^3}\right)$. By equation \eqref{eq:l_in_I} the probability that $l \notin I$ is $o\left(\frac{1}{n^3}\right)$. Combining these observations we find that for all $\epsilon >0$ the terms in \eqref{eq:petrov} are vanishing:
\begin{align}
\label{eq:njk_exp_eps}
\p{\left\lvert \Np{j}{k} - \mathbb{E}\left[\Np{j}{k}\right] \right\rvert > n\epsilon} = o\left(\frac{1}{n^3}\right),
\end{align}
  which in turn proves that the limit in equation \eqref{eq:con_dis_bond} holds $\measure$-a.s.
  
   It remains to show that the first, first mixed and second moments of $\frac{\Np{j}{k}}{n}$ converge $\measure$-a.s. to those of $\diprobdisbond$.      In Section~\ref{subsubsc:degree_dis_bond} we showed that $\sum_{j,k=0}^\infty j \pbond{j}{k}= \sum_{j,k=0}^\infty k \pbond{j}{k}$, and due to the graph context of the problem, $\sum_{j,k=0}^\infty j \Np{j}{k} = \sum_{j,k=0}^\infty k \Np{j}{k}$. 
Therefore, it is sufficient to show that one of the first moments converges. Let additionally $\Qn:= \frac{1}{n}\sum_{j,k=0}^\infty j \Np{j}{k}$, we will then show that
$
Q:=\lim_{n \rightarrow \infty} \Qn=\per\mum,\; \measure\text{-a.s.}
$
Let
$ \Xkn := \frac{1}{n}\sum_{j=0}^{\kappa}\sum_{k=0}^{\kappa} j \Np{j}{k}
$ and remark that $\Xkn \leq \Qn$. 
Since $\diarray$ is a \proper degree progression, for all $\epsilon >0$  there exists $\widetilde{\kappa}\left(\epsilon\right), \tilde m\left(\epsilon\right)$ such that for all   $ \kappa > \widetilde{\kappa}$ and $n > \tilde m$
\begin{align}
\label{eq:n_kappa_epsilon}
\frac{1}{n} \sum_{\mathclap{\substack{(d^-,d^+) = (0,0),\\
			d^- > \kappa  \,\text{or}\, d^+ > \kappa}}}^{(d_{\max}, d_{\max})} j\N{j}{k} < \epsilon.
\end{align}
Thus for $\kappa > \widetilde{\kappa}$ there holds  $ \Xkn \leq  \Qn \leq \Xkn + \epsilon$. This implies that if
\begin{align}
\label{eq:lim_X_bond} 
\Xk:=\lim_{n\rightarrow \infty } \Xkn = \sum_{j=0}^{\kappa}\sum_{k=0}^{\kappa} j\pbond{j}{k} , \; \text{ $\measure$-a.s.},
\end{align}
then $\lim\limits_{n\rightarrow \infty } \Qn = \Q$ holds $\measure$-a.s. as well~\cite[$(3.5)$]{fountoulakis2007percolation}. Thus the goal is to prove equation \eqref{eq:lim_X_bond}. Applying \cite[Lemma $6.8$]{petrov1995limit} we can show that this limit holds $\measure$-a.s., if for any $\epsilon >0$, 
\begin{align}
\label{eq:sum_prob_finite}
\sum_{n=1}^\infty \p{\left\lvert \Xkn - \Xk \right \rvert > \epsilon} < \infty.
\end{align}
We show this analogously to the proof of equation \eqref{eq:petrov}. By the definition of $
\Xkn, \Xk$ and $\pbond{j}{k}$, for all $\epsilon >0$ there exists $\widetilde{N}$ such that for all $n >\widetilde{N}$,
$
\left\lvert \mathbb{E}\left[\Xkn\right] - \Xk \right\rvert < \frac{\epsilon}{2}. 
$
Combing this with the reverse triangle inequality, we find for any $\epsilon >0$,
\begin{align*}
\p{\left\lvert \Xkn - \Xk \right \rvert > \epsilon} \leq \p{\left\lvert \Xkn - \mathbb{E}\left[\Xkn\right] \right \rvert > \frac{\epsilon}{2}}.
\end{align*}
Remark that 
$
\left\lvert \Xkn - \mathbb{E}\left[\Xkn\right] \right \rvert = \frac{1}{n}\sum_{j=0}^\kappa \sum_{k=0}^\kappa j \left(\Np{j}{k} - \mathbb{E}\left[\Np{j}{k}\right] \right).
$
This implies that for  $\epsilon' = \frac{\epsilon}{2\sum_{j \leq \kappa}j}$ there holds
\begin{align*}
\p{\left\lvert \Xkn - \mathbb{E}\left[\Xkn\right] \right \rvert > \frac{\epsilon}{2}} \leq \sum_{j\leq \kappa, k \leq \kappa} \p{\frac{1}{n}\lvert \Np{j}{k} - \mathbb{E}\left[\Np{j}{k}\right] \rvert > \epsilon '}.
\end{align*}
Using  equation \eqref{eq:njk_exp_eps} we find
$
\p{\left\lvert \Xkn - \mathbb{E}\left[\Xkn\right] \right \rvert > \frac{\epsilon}{2}} \leq \sum_{j\leq \kappa, k \leq \kappa} o\left(\frac{1}{n^3}\right) \leq o\left(\frac{1}{n^{2\frac{7}{9}}}\right).
$
Here we used the fact that $d_{\max} = O\left(n^{1/9}\right)$ and that for $j >d_{\max}$ or $k>d_{\max}$ or both, there holds $\Np{j}{k} = \mathbb{E}\left[\Np{j}{k}\right]=0$. 
This proves  \eqref{eq:sum_prob_finite} and hence proves that $\Qn$ converges $\measure$-a.s. to $\Q$.
 By redefining $\Qn, \Q, \Xkn, \Xk$ and setting 
$
 \epsilon': = \min \left \{ \frac{\epsilon}{2}, \frac{\epsilon}{2\sum_{j \leq \kappa}j}, \frac{\epsilon}{2\sum_{j \leq \kappa, k \leq \kappa }jk}, \frac{\epsilon}{2\sum_{j \leq \kappa}j^2}, \frac{\epsilon}{2\sum_{k \leq \kappa}k^2} \right \},
$
 similar derivations hold for the first mixed moment and the second moments, that is one may show that all the moments of interest and the distribution itself converge simultaneously $\measure$-a.s. for an element of $\D$. Thus we have shown that $\pdiarray$  is $\measure$-a.s. feasible.
  \end{proof}
  \end{Lemma}
  
  \begin{Lemma}\label{lm:bondproper}
The degree progression after bond percolation is proper $\measure$-a.s.
\begin{proof}
The degree progression $\pdiarray$  is feasible due to Lemma~\ref{lm:bondfeasible}, thus to
  prove that $\pdiarray$   is $\measure$-a.s. proper according to Definition~\ref{def:proper_deg_progression}, it remains to show that for a feasible $ \pdiarray \in \D$:
  the maximum degree  of the percolated degree sequence is bounded $d^{\per}_{\max} \leq \frac{n^{1/12}}{\ln n}$
   and $$\rho^{\per}(n):= \max \left\{ \frac{\sum_{j,k=0}^\infty j^2k\Np{j}{k}}{\mum^{\per} n},  \frac{\sum_{j,k=0}^\infty jk^2\Np{j}{k}}{\mum^{\per} n}\right\} = o \left(\frac{n^{1/12}}{\ln n}\right).$$
        As the degree progression before percolation $\diarray$ is proper, there holds  $d_{\max} \leq \frac{n^{1/12}}{\ln n}$. Percolation can only decrease the in-degree  and the out-degree of a vertex, implying $d_{\max}^{\per} \leq d_{\max}$.  Together these observations show that $d^{\per}_{\max} \leq \frac{n^{1/12}}{\ln n}$ for any  $\pdiarray \in D$.
It remains to show $\rho^{\per} (n)= o \left(\frac{n^{1/12}}{\ln n}\right)$. Consider the sum $\sum_{j,k=0}^\infty j^2k\Np{j}{k}$ (respectively $ \sum_{j,k=0}^\infty jk^2\Np{j}{k}$). Each vertex of degree $(j,k)$ contributes $j^2k$ (or $jk^2$) to the total sum. As the in-degree and the out-degree can only decrease due to percolation, this implies that $$\sum_{j,k=0}^\infty j^2k\Np{j}{k} \leq \sum_{j,k=0}^\infty j^2k\N{j}{k} \quad \text{and} \quad \sum_{j,k=0}^\infty jk^2\Np{j}{k} \leq \sum_{j,k=0}^\infty jk^2\N{j}{k}.$$
Recall that $\mum^{\per, \text{bond}} = \per \mum$ and $\per$ is a constant.  These observations, together with the fact that $\diarray$ is \proper\!\!, imply  that
$
  \rho^{\per}(n) \leq
  \frac{\rho(n)}{\per} = o\left(\frac{n^{1/12}}{\ln n}\right).
$
  Thus $\pdiarray \in D$ is $\measure$-a.s. proper.
  \end{proof}
    \end{Lemma}

  Let $E \subset \D$ be the event over which the degree progression is proper. As any element of $\D$ is $\measure$-a.s. proper, there holds  $\measure\left(E\right) =1$. 
  For any $\pdiarray \in E$ we may apply  Theorem~\ref{thm:giant_scc} to a sequence of random multigraphs $\mpgrapharray$ arising from uniformly random configurations. Recall that Lemma~\ref{lm:percolation_is_CM} implies that this is the case for all $n$ if we condition on  $\diseqrv=\diseqp$. We will now fix $\pdiarray \in E$ and apply Theorem ~\ref{thm:giant_scc} to $\mpgrapharray$, distinguishing two cases for $\per$: $ \per < \hat{\per}^{\text{bond}}$ and $ \per > \hat{\per}^{\text{bond}}$, with  $\hat{\per}^{\text{bond}} = \frac{\mum}{\mum[11]}$.

Let $ \per < \hat{\per}^{\text{bond}}$. Define $\mathcal{A}_\epsilon\left(\diseqp\right)$ to be the set of all multigraphs  obeying $\diseqp$ for which the largest strongly connected component contains no more than $\epsilon_1 n$ vertices for any  $\epsilon_1 \in (0,1)$. As $\frac{\mum[11]^{\per}}{\mum^{\per}} = \per\frac{\mum[11]}{\mum} <1$, Theorem~\ref{thm:giant_scc} states that for all $\epsilon_1$:
  \begin{align}
  \label{eq:giant_not_bond}
  \lim_{n \rightarrow \infty} \p{\mpgraph \in \mathcal{A}_{\epsilon_1}\left(\diseqp\right) \,\vert\,  \diseqrv=\diseqp} = 1.
  \end{align}
  Next consider $ \per > \hat{\per}^{\text{bond}}$.  Define $\mathcal{B}_{\epsilon}\left(\diseqp\right)$ to be the set of all graphs whose largest strongly connected component contains $\epsilon_2 n$ vertices, with $\epsilon_2 \in (0,1)$. Since  $\frac{\mum[11]^{\per}}{\mum^{\per}} = \per\frac{\mum[11]}{\mum} >1$, Theorem~\ref{thm:giant_scc} states that there exists a unique $\epsilon_2=\csb$ such that the value of the following limit is nonzero:
  \begin{align}
    \label{eq:giant_exits_bond}
  \lim_{n \rightarrow \infty} \p{\mpgraph \in \mathcal{B}_{\epsilon_2}\left(\diseqp\right) \,\vert\, \diseqrv = \diseqp } = 1.
  \end{align}
  Moreover, the theorem determines that $\csb :=1-U_{\per}^\text{bond}(x^*,1)-U_{\per}^\text{bond}(1,y^*)+U_{\per}^\text{bond}(x^*,y^*),$
where $x^*,y^*\in(0,1)$ are the unique solutions of
\begin{equation}
\begin{aligned}
\label{eq:xystar}
x^*= U^-_{\per}(x^*),\\
y^*= U^+_{\per}(y^*),
\end{aligned}
\end{equation}
and functions $U_{\per}^+,U^-_{\per}$ and $U_{\per}^\text{bond}$ are as defined in \eqref{eq:gfs}.

 To finalise the proof for Theorems~\ref{thm:percolation_di_strong} and~\ref{thm:percolation_di_strong2}, we need to supplement Equations  \eqref{eq:giant_exits_bond} and  \eqref{eq:giant_not_bond} with two  minor observations. First,  the theorem is stated for a percolated multigraph progression $\mpgrapharray$ without conditioning on the degree progression of the percolated graphs. As $\measure\left(E\right) =1$, an analogous argument to that of Fountoulakis~\cite[p. $348$]{fountoulakis2007percolation}  is applied to show that:
\begin{itemize}
	\item   If $\per < \hat{\per}^{\text{bond}}, \;
	\lim \limits_{n \rightarrow \infty} \p{\mpgraph \in \mathcal{A}_\epsilon\left(\diseqp\right)} = 1$ for all $\epsilon \in (0,1)$.
	\item   If $\per > \hat{\per}^{\text{bond}},\;
	\lim\limits_{n \rightarrow \infty} \p{\mpgraph \in \mathcal{B}_{\csb}\left(\diseqp\right)} = 1$\\ and  $
	\lim\limits_{n \rightarrow \infty} \p{\mpgraph \in \mathcal{B}_{\epsilon}\left(\diseqp\right)} = 0$ for all $\epsilon \in (0,1), \epsilon \neq \csb$ .
\end{itemize}
Second,  Theorems~\ref{thm:percolation_di_strong} and~\ref{thm:percolation_di_strong2} make assertions about uniformly random \emph{simple} graphs instead of random multigraphs. Replace the
the graph $\mgraphseq$ in  Lemma~\ref{lm:multi_simple_0}  by the graph $\mpgraph$ and condition on the graph to which percolation is applied ($\graphseq$) being simple. This yields a slightly different variant of the lemma that does not require additional changes to the proof. Now applying this variant of Lemma~\ref{lm:multi_simple_0}  to the above limits, we deduce an equivalent statement as above for $\pgraph$:
\begin{itemize}
	\item  If $\per < \hat{\per}^{\text{bond}}, \;
	\lim\limits_{n \rightarrow \infty} \p{\pgraph \in \mathcal{A}_\epsilon\left(\diseqp\right)} = 1$ for all $\epsilon \in (0,1)$.
	\item  If $\per > \hat{\per}^{\text{bond}},$ $
	\lim\limits_{n \rightarrow \infty} \p{\pgraph \in \mathcal{B}_{\csb}\left(\diseqp\right)} = 1$\\ and $
	\lim\limits_{n \rightarrow \infty} \p{\pgraph \in \mathcal{B}_{\epsilon}\left(\diseqp\right)} = 0$ for all $\epsilon \in (0,1), \epsilon \neq \csb$.
\end{itemize}
completing the proofs of Theorems~\ref{thm:percolation_di_strong} and~\ref{thm:percolation_di_strong2} for the case of bond percolation. 
\subsection{Site percolation}
\label{subsc:site}
The proof of Theorems~\ref{thm:percolation_di_strong} and~\ref{thm:percolation_di_strong2} for site percolation has a similar structure as that for bond percolation.
  Hence, we will refer back to Section~\ref{subsc:bond} where applicable.  The proof is split into three steps. First, in Section~\ref{subsubsc:random_site} we show that applying site percolation to a uniformly random  configuration results in another uniformly random  configuration if we condition on the degree sequence after percolation.
Second, we determine the limit of the expected number of vertices with degree $(j,k)$ after site percolation, see Section~\ref{subsubsc:degree_dis_site}. 
The proof is completed in Section~\ref{subsubsc:theorem_site} by combining the first two steps with the results of Section~\ref{subsc:bond}.

 \noindent Recall from Section~\ref{sc:main} that deleting a vertex means that we remove all edges adjacent to this vertex. In the setting of the configuration model this implies that all stubs attached to a deleted vertex are removed. Let us denote these stubs by  $\left(\instubr, \outstubr\right)$. As site percolation removes any edges adjacent to a vertex, the match of any stub in $\left(\instubr, \outstubr\right)$ will be removed as well. A stub in  $\left(\instubr, \outstubr\right)$ may or may not have its match in the same set, as it might happen that both endpoints of one edge are deleted. Let 
$\left(\instubm, \outstubm\right)$ contain all the matches of stubs in $\left(\instubr, \outstubr\right)$ that are not connected to a deleted vertex. Thus $\instubr \cup \instubm$ (respectively $\outstubr \cup \outstubm$) are all in-stubs (out-stubs) removed by site percolation. 
The stubs that survive percolation are still denoted by $\left(\instubp, \outstubp\right)$. Remark that this implies
\begin{align}
\label{eq:site_union}
\instub = \instubp \cup \instubr\cup \instubm \quad \text{and} \quad  \outstub = \outstubp \cup \outstubr\cup \outstubm.
\end{align} These definitions of $\left(\instubr, \outstubr\right)$ and $\left(\instubm, \outstubm\right)$ will be used throughout the proof.

\subsubsection{A site-percolated configuration is a uniformly random configuration}
\label{subsubsc:random_site}
We show in Lemma~\ref{lm:site_cm_points_percolation} that conditional on the stubs that are removed by site percolation, the matching on the surviving stubs is uniformly random.  This lemma, in turn, allows us to formulate Lemma~\ref{lm:site_cm} stating that conditional on the degree sequence after percolation, there remains a uniformly random configuration.
\begin{Lemma}
	\label{lm:site_cm_points_percolation}
	Apply site percolation to a uniformly random configuration $\M$ on $\left(\instub, \outstub\right)$. Conditional on the elements of $\left(\instubr, \outstubr\right)$ and $\left(\instubm,\outstubm\right)$, each configuration on $\left(\instubp,\outstubp\right)$ is equally likely.
	\begin{proof}
		According to equation \eqref{eq:site_union}, fixing the elements of $\left(\instubr, \outstubr\right)$ and $\left(\instubm,\outstubm\right)$, uniquely determines the elements of $\left(\instubp,\outstubp\right)$. Choosing the elements of $\left(\instubr, \outstubr\right)$ and $\left(\instubm,\outstubm\right)$ furthermore implies that the configuration $\M$ is the union of a configuration on $\left(\instubr\cup \instubm, \outstubr\cup \instubm\right)$  with the one on $\left(\instubp,\outstubp\right)$. As $\M$ is a uniformly random configuration obeying this split and the elements of $\left(\instubp,\outstubp\right)$ are fixed, the configuration on $\left(\instubp,\outstubp\right)$ will be a uniformly random one.

	\end{proof} 
\end{Lemma}
\begin{Lemma}
	\label{lm:site_cm}
	Apply site percolation to a uniformly random configuration  on
	$\left(\instub, \outstub\right)$.  Conditional on $\diseqrv = \diseqp$, any configuration on $\left(\instub[\diseqp], \outstub[\diseqp]\right)$ is equally likely. 
	\begin{proof}
		Define $l = \lvert \instubp \rvert$ and let $S(\diseqp)$ contain all sets of surviving stubs $\left(\instubp,\outstubp\right)$ that induce the degrees sequence $\diseqp$. Fix a matching $\Mp$ of $\left(\instub[\diseqp], \outstub[\diseqp]\right)$. Then it holds that
		\begin{align*}
		\p{\Mp \middle| \diseqrv = \diseqp} = \sum_{\left(A,B\right)\in S(\diseqp)} &	\p{\Mp \middle| \diseqrv =
	 \diseqp,  \left(\instubp, \outstubp\right) =
	  \left(A,B\right) } \times \\
	  &\p{ \left(\instubp, \outstubp\right) = \left(A,B\right)| \diseqrv = \diseqp }.
		\end{align*}
		Remark that $\p{\Mp \middle| \diseqrv = \diseqp,  \left(\instubp, \outstubp\right) = \left(A,B\right) } = \p{\Mp \middle| \left(\instubp, \outstubp\right) = \left(A,B\right) }$  as $(A,B) \in S(\diseqp)$ implies that $\left(\instubp, \outstubp\right)$ induces degree sequence $\diseqp$. 
		Using Lemma~\ref{lm:site_cm_points_percolation} and the bijection between $\left(\instubp, \outstubp \right)$ and $\left(\instub[\diseqp], \outstub[\diseqp]\right)$ we find:
		   $$\p{\Mp \middle|  \left(\instubp, \outstubp\right) = \left(A,B\right) } = \frac{1}{l!} .$$
		   Furthermore, combining these observations with $$\sum_{\left(A,B\right)\in S(\diseqp)}  \mathbb{P}\left[\left(\instubp, \outstubp\right)= \left(A,B\right) | \diseqrv= \diseqp \right] = 1$$ following from the definition of $S(\diseqp)$, we obtain
		\begin{align*}
		\p{\Mp \middle| \diseqrv = \diseqp}	&= \frac{1}{l!} \sum_{\left(A,B\right)\in S(\diseqp)}  \mathbb{P}\left[  \left(\instubp, \outstubp\right) = \left(A,B\right)| \diseqrv= \diseqp \right] = \frac{1}{l!},
		\end{align*} 
		completing the proof.
	\end{proof} 
\end{Lemma}

\subsubsection{The expected number of vertices with degree $(j,k)$ after site percolation}
\label{subsubsc:degree_dis_site}
The next step in the proof of Theorem~\ref{thm:percolation_di_strong} for site percolation is proving existence of the limit of $\frac{\Np{j}{k}}{n}$
 and determining its value, $\psite{j}{k}$.   Then, given $\psite{j}{k}$,  $\hat{\per}^{\text{site}}$ is determined and, in Section~\ref{subsubsc:theorem_site}, $\hat{\per}^{\text{site}}$ is shown  to be the desired threshold for site percolation. 
 
\begin{Lemma}
	\label{lm:Nkj_pi_site}
	Let $\Np{j}{k}$ be the number of vertices in the site-percolated graph (or configuration) with in-degree $j$ and out-degree $k$. The following limit exists
\begin{align}
\label{eq:lim_p_per_site}
\psite{j}{k}:=\lim_{n\rightarrow \infty} \frac{\mathbb{E}\left[\Np{j}{k}\right]}{n}, \text{ for } j,k\in\mathbb N_0,
\end{align}
and
\begin{align}
\label{eq:psite_pbond}
\psite{j}{k} = \begin{cases}
\per\pbond{j}{k},&(j,k)\neq(0,0),\\
\per\pbond{0}{0} + 1-\pi,& (j,k)=(0,0),
\end{cases}
\end{align}
where $\pbond{j}{k}$ is defined  in Lemma~\ref{lm:Nkj_pi}.
\begin{proof}
If  $j, k>d_{\max}$ then $\psite{j}{k}=\Np{j}{k}=\N{j}{k}=0$.
Consider $0\leq j,k \leq d_{\max}$. 
Let $\Nr{d^-}{d^+}$ be the number of vertices of degree $(d^-,d^+)$ before percolation that are not deleted.  Thus $\N{d^-}{d^+}- \Nr{d^-}{d^+}$ is the number of deleted vertices of degree $(d^-,d^+)$. Each vertex is deleted with probability $1-\per$ independently of other vertices, hence:
\begin{align}
\label{eq:expect_Nr}
&\mathbb{E}\left[\Nr{d^-}{d^+}\right] = \per\N{d^-}{d^+},\\
&\mathbb{E}\left[\N{d^-}{d^+}- \Nr{d^-}{d^+}\right] =  \left(1-\per\right)\N{d^-}{d^+}.
\end{align}
 A deleted vertex will have degree $(0,0)$ after percolation with probability $1$. Let $P_{j,k}\left(d^-,d^+\right)$ be the probability that a non-deleted vertex of degree $(d^-,d^+)$  has degree $(j,k)$ after percolation. 
For $(j,k) = (0,0)$ we have
\begin{align}
\label{eq:def_exp_Dp_00}
\mathbb{E}\left[\Np{0}{0}\right] = \sum_{d^- = 0}^{d_{\max}}\sum_{d^+ = 0}^{d_{\max}}\left( \left(1-\per\right) \N{d^-}{d^+} + \per P_{0,0}\left(d^-,d^+\right) \N{d^-}{d^+}\right),
\end{align}
and otherwise, 
\begin{align}
\label{eq:def_exp_Dp}
\mathbb{E}\left[\Np{j}{k}\right] = \sum_{d^- =j}^{d_{\max}}\sum_{d^+ =k}^{d_{\max}} \per P_{j,k}\left(d^-,d^+\right) \N{d^-}{d^+}.
\end{align}
Hence we need to derive the expression for $P_{j,k}\left(d^-,d^+\right)$. 
Let $$\si: = \lvert \instubp \cup \instubm \rvert,\; \so:  = \lvert\outstubp \cup \outstubm \rvert,\;  \ri := \lvert \instubm\rvert,\; \ro:= \lvert \outstubm \rvert.$$ Note that there must hold $\si - \ri = \so - \ro $ as $\si - \ri  = \lvert \instubp \rvert $, $\so - \ro  = \lvert \outstubp \rvert $ and the remaining configuration on $\left(\instubp, \outstubp\right)$ forms a directed graph. Let $P_{j,k}\left(d^-,d^+, \si, \so, \ri, \ro\right)$ denote the probability  $P_{j,k}\left(d^-,d^+\right)$ conditional on the values $\si,\so,\ri,\ro$. We will now determine this conditional probability.  Site percolation combines the independent random processes of deleting vertices and creating a uniformly random configuration on $\left(\instub, \outstub\right)$. As these processes are independent, we may first determine the elements of $\left(\instubr,\outstubr\right)$ and then randomly create a configuration on $\left(\instub, \outstub\right)$. Thus conditional on  the value $\ri$ (respectively $\ro$), each subset of $\instub \setminus \instubr$($\outstub \setminus \outstubr$)  of this size is equally likely to be $\instubm$ (or $\outstubm$). This implies that
\begin{align}
\label{eq:prob_remain_initial}
P_{j,k}\left(d^-,d^+, \ri, \ro, \si, \so\right) = {{d^-}\choose {d^- -j}} {{d^+}\choose {d^+ -k}} \frac{{{\si - d^-}\choose {\ri - d^- +j}}}{{{\si}\choose{\ri}}}\frac{{{\so - d^+}\choose {\ro - d^+ +k}}}{{{\so}\choose{\ro}}}.
\end{align} 
 To approximate this probability we formulate the following propositions that show that with high probability $\si,\so$  remains in some bounded interval $I'_n$ and $\ri,\ro$ in $I_n$. 

 \begin{Proposition}\label{prop:I'}
 Let $I'_n := \left[m\per - n^{2/3}\ln n, m\per + n^{2/3}\ln n\right] $.
	There holds 		
	$$\p{ \si,\so \in I'_n }=1-e^{-\Omega\left(\ln^2 n\right)}$$
	 for both $\si$ and $\so$  separately.
 \begin{proof}
 By using equation \eqref{eq:expect_Nr}, we obtain:
$
\mathbb{E}\left[\si\right] = \sum_{d^-=0}^{d_{\max}}\sum_{d^+=0}^{d_{\max}} \per d^-\N{d^-}{d^+} = m\per$ and $\mathbb{E}\left[\so\right] = \sum_{d^-=0}^{d_{\max}}\sum_{d^+=0}^{d_{\max}} \per d^+ \N{d^-}{d^+} =m\per .
$
Using $d_{\max} \leq n^{1/9}$ and Hoeffding's inequality we also find that:
\begin{align}
\label{eq:bound_s}
\p{\left\lvert \si - \mathbb{E}\left[\si\right]\right\rvert > n^{2/3}\ln n} \leq e^{-\Omega\left( \ln^2 n\right)} \; \text{and}\; \p{\left\lvert \so - \mathbb{E}\left[\so\right]\right\rvert > n^{2/3}\ln n} \leq e^{-\Omega\left( \ln^2 n\right)},
\end{align}
proving the claim.
\end{proof}
 \end{Proposition}
\begin{Proposition}
	\label{lm:size_b-_b+}
	Let $I_n := \left[m\per(1-\per) - n^{2/3}\ln^2 n, m\per(1-\per) + n^{2/3}\ln n^2\right]$.
	There holds 
	$$\p{ \ro,\ri \in I_n \mid \si, \so \in I'_n}=1-e^{-\Omega\left(\ln^2 n\right)}$$
	 for both $\ro$ and $\ri$  separately.
	 	\begin{proof}
		We present the proof for $\ri$. The proof for $\ro$ is identical to the one for $\ri$ up to switching the roles of in-stubs and out-stubs. 
		Since we consider a uniformly random configuration on $\left(\instub, \outstub\right)$, the probability that any in-stub  is matched to  an out-stub in $\outstubr$ is $\frac{m-\so}{m}= (1-\per) \left( 1+ \mathcal{O}\left(n^{-1/3}\ln n\right)\right)$  as $\si,\so \in I'_n$. Since $\ri$ equals the number of in-stubs in $\instub \setminus\instubr$ with a match in $\outstubr$, this implies that
		\begin{align*}
		\mathbb{E}\left[\ri\right] = \si \frac{m-\so}{m} = m\per\left(1-\per\right) \left(1 + \mathcal{O} \left(n^{2/3}\ln n\right)\right). 
		\end{align*} 
		To complete the proof, we will now show that
		\begin{align*}
		\p{\lvert  \ri - \mathbb{E}\left[\ri\right] \rvert  > n^{2/3}\ln^2 n } \leq e^{- \Omega\left(\ln^2 n\right)}.
		\end{align*}
		This is realised by applying Theorem~\ref{thm:mcdiarmid} to the space of configurations on $\left(\instub, \outstub \right)$ with the symmetric difference as the metric. The value of $\ri$ plays the role of the function $f$.  To partition this space, we order the in-stubs of $\instub$. Define an $i$-prefix to be the first $i$ in-stubs together with their match. An element of the partition $\mathcal{P}_k$ consists of all  configurations with the same $k$-prefix for all $k \in \{0,1,\ldots, m\}$. For any $A,B \in \mathcal{P}_k$ such that $A,B\subset C \in \mathcal{P}_{k-1}$ a bijection $\phi: A \rightarrow B$ can be defined. Denote the $k^{\text{th}}$ pair of a configuration in $A$ by $(x,y_A)$ and the $k^{\text{th}}$ pair of a configuration in $B$ by $(x,y_B)$. Then $\phi$ maps $\M \in A$ to the configuration in $B$ with $(x,y_A)$ replaced by $(x,y_B)$ and with $y_A$ the match of the in-stub in $\M$ matched to $y_B$. By definition of $\phi$ it follows that
		 $c_k:=\lvert \M -\phi(\M)\rvert = 4 $ for all $k \in \{1,2,\ldots ,m\}$. As the value of $\ri$  changes by at most the symmetric difference of the two matchings, Theorem~\ref{thm:mcdiarmid} states that:
		\begin{align*}
		\p{\lvert  \ri - \mathbb{E}\left[\ri\right] \rvert  > n^{2/3}\ln^2 n } \leq  2\exp \left(\frac{n^{4/3}\ln^2 n}{2m}\right) = e^{- \Omega\left(\ln^2 n\right)},
		\end{align*}
		as $m \leq nd_{\max} \leq n^{10/9}$. 
	\end{proof}
\end{Proposition}	
  For $\ri,\ro \in I_n$ and $ \si,\so \in I'_n$  there holds uniformly:
\begin{align*}
{d \choose {d-i}}\frac{{{s - d} \choose {r -d +i}}}{{{s}\choose {r}}} = {d \choose {d-i}} \left(1-\per\right)^{d-i}\per^i\left(1 + \mathcal{O}\left(\frac{\ln^2 n}{n^{1/3}}\right)\right).	\end{align*}
Plugging this equation into \eqref{eq:prob_remain_initial} yields that 
\begin{align*}
P_{j,k}\left(d^-,d^+, \ri, \ro, \si, \so\right) = {{d^-}\choose {d^- -j}} {{d^+}\choose {d^+ -k}} \per^{j+k}\left(1-\per\right)^{d^-+d^+-j-k} \left(1 + \mathcal{O}\left(\frac{\ln^2 n}{n^{1/3}}\right)\right),
\end{align*}
holds uniformly for all $\si, \so \in I'_n$ and $\ri,\ro \in I_n$.
However we cannot yet determine this probability if at least one of the following conditions is violated: $\si, \so \in I'_n$, or $\ri,\ro \in I_n.$
 We show that such violations are unlikely, that is, instead of bounding the probability 
$\p{\si \notin I'_n \,\vee	\, \so \notin I'_n \,\vee	\,\ri \notin I_n \,\vee	\,\ro \notin I_n},$ we add a condition on  $\Nr{d^-}{d^+}$,  allowing us to bound the value of $\mathbb{E}\left[\Np{j}{k}\right]$. Theorem~\ref{thm:hoeffding} implies that
\begin{align}
\label{eq:bound_Nr}
\p{\Nr{d^-}{d^+} - \mathbb{E}\left[\Nr{d^+}{d^-}\right] \rvert > \sqrt{n}\ln n} < e^{-\Omega\left(\ln^2 n\right)}.
\end{align}
In combination with equation \eqref{eq:expect_Nr} this implies that  
$$\Nr{d^-}{d^+} \in I''_n(d^-,d^+) := \left[\max\left\{\per\N{d^-}{d^+} - \sqrt{n}\ln n, 0\right\}, \per\N{d^-}{d^+} + \sqrt{n}\ln n\right],$$ with probability $1-e^{-\Omega\left(\ln^2 n\right)}$. Combining this with Propositions~\ref{prop:I'}~and~\ref{lm:size_b-_b+} gives:
\begin{align*}
&\!\!\!\!\!\!\p{\si \notin I'_n \,\lor\, \so \notin I'_n \,\lor\,\ri \notin I_n \,\lor \notin I_n  \lor \Nr{d^-}{d^+} \notin I''_n\left(d^-,d^+\right)} \;\;\\
\leq & \p{\si \notin I'_n} + \p {\so \notin I'_n} + \p{\ri \notin I_n} + \p{\ro \notin I_n}  + \p{\Nr{d^-}{d^+} \notin I''_n\left(d^-,d^+\right)}\\
= & o\left(\frac{1}{n^3}\right) + \p{\ri \notin I_n} + \p{\ro \notin I_n}.
\end{align*}
By the law of total probability:  
\begin{equation*}
\begin{split}
\p{\ri \notin I_n} =&\p{\ri \notin I_n \mid \si \in I'_n, \so \in I'_n}\p{ \si \in I'_n, \so \in I'_n}+ \\
& \p{\ri \notin I_n \mid \si \notin I'_n, \so \in I'_n}\p{ \si \notin I'_n, \so \in I'_n} +\\
& \p{\ri \notin I_n \mid \si \in I'_n, \so \notin I'_n}\p{ \si \in I'_n, \so \notin I'_n}+\\
 & \p{\ri \notin I_n \mid \si \notin I'_n, \so \notin I'_n}\p{ \si \notin I'_n, \so \notin I'_n} = o\left(n^{-3}\right).
\end{split}
\end{equation*}
In a similar fashion, it is shown that $\p{\ro \notin I_n} =  o\left(n^{-3}\right)$.
Thus there holds
\begin{align}
\label{eq:prob_out_I}
\p{\si \notin I'_n \lor \so \notin I'_n \lor \ri \notin I_n \lor\ro \notin I_n  \lor \Nr{d^-}{d^+} \notin I''_n\left(d^-,d^+\right)} =  o\left(n^{-3}\right). 
\end{align}
This allows us to determine a lower and upper bound for the value $\mathbb{E}\left[\Np{j}{k}\right]$. Since $\Nr{d^-}{d^+}  \leq \N{d^-}{d^+}$ and  $\diarray$ is proper, for all $\epsilon > 0$ there exist $\kappa\left(\epsilon\right)$ and $N\left(\epsilon\right)$ such that for all $n > N:$
\begin{align}
\label{eq:epislon_nd}
\sum_{\mathclap{\substack{(d^-,d^+) = (0,0)\\
			d^- \geq \kappa +1 \,\text{or}\, d^+ \geq \kappa +1}}}^{(d_{\max}, d_{\max})}
		P_{j,k}\left(d^-,d^+\right) \Nr{d^-}{d^+} \leq  \sum_{\mathclap{\substack{(d^-,d^+) = (0,0)\\
					d^- \geq \kappa +1 \,\text{or}\, d^+ \geq \kappa +1}}}^{(d_{\max}, d_{\max})}
				 \N{d^-}{d^+} \leq \epsilon n.
\end{align}
In combination with equation \eqref{eq:def_exp_Dp} this implies for $(j,k) \neq (0,0)$
\begin{align}
\label{eq:bound_Dp_first}
\sum_{ \substack{d^- =j,\\ d^+ =k}}^{\kappa}P_{j,k}\left(d^-,d^+\right) \Nr{d^-}{d^+} \leq \mathbb{E}\left[\Np{j}{k} \right]\leq  \sum_{ \substack{d^- =j,\\ d^+ =k}}^{\kappa}P_{j,k}\left(d^-,d^+\right) \Nr{d^-}{d^+} + \epsilon n.
\end{align} 
Using equation \eqref{eq:prob_out_I} on the left-hand side of the above equation we find
\begin{equation*}
\begin{split}
\mathbb{E}\left[\Np{j}{k}\right] \geq& \sum_{\substack{d^-=j,\\d^+=k}}^{\kappa}\; \sum_{\substack{\tilde r^- \in I'_n,\\\tilde r^+\in I'_n}} \;\sum_{\substack{\tilde s^-\in I_n,\\\tilde s^+ \in I_n}}\;\sum_{\tilde{d}_{d^-,d^+} \in I''_n(d^-,d^+)} \tilde{d}_{d^-,d^+}P_{j,k}\left(d^-,d^+,\tilde r^-,\tilde r^+,\tilde s^-,\tilde s^+\right)\\
&\times\p{\ri=\tilde r^-, \ro=\tilde r^+, \si = \tilde s^-, \so = \tilde s^+, \Nr{d^-}{d^+} = \tilde{d}_{d^-,d^+}} + o \left(\frac{1}{n^2}\right).
\end{split}
\end{equation*}
As equation \eqref{eq:bound_Nr} implies that 
$$\sum_{\tilde{d}_{d^-,d^+} \in I''_n(d^-,d^+)}\tilde{d}_{d^-,d^+}\p{\Nr{d^-}{d^+} = \tilde{d}_{d^-,d^+}} =  \mathbb{E}\left[\Nr{d^-}{d^+}\right] + o\left(\frac{1}{n^2}\right),$$
 we obtain the lower bound:
\begin{align*}
\mathbb{E}\left[\Np{j}{k}\right] \geq o\left(\frac{1}{n^2}\right)+ \per \sum_{d^-=j}^{\kappa}\sum_{d^+=k}^{\kappa} \N{d^-}{d^+} {{d^-}\choose {d^- -j}} {{d^+}\choose {d^+ -k}}  \times\\
\per^{j+k}\left(1-\per\right)^{d^-+d^+-j-k}
\left(1 + \mathcal{O}\left(\frac{\ln^2 n}{n^{1/3}}\right)\right) .
\end{align*}
In a similar fashion, using the right-hand side of equation \eqref{eq:bound_Dp_first}, gives the upper bound:
\begin{align*}
\mathbb{E}\left[\Np{j}{k}\right] \leq \epsilon n +  o\left(\frac{1}{n^2}\right)+\per \sum_{d^-=j}^{\kappa}\sum_{d^+=k}^{\kappa} \N{d^-}{d^+} {{d^-}\choose {d^- -j}} {{d^+}\choose {d^+ -k}}
 \times\\
  \per^{j+k}\left(1-\per\right)^{d^-+d^+-j-k}
\left(1 + \mathcal{O}\left(\frac{\ln^2 n}{n^{1/3}}\right)\right) .
\end{align*}
Combining the upper and lower bounds together proves  convergence of the limit for $j,k>0$:
\begin{align}
\label{eq:psite}
\psite{j}{k}:=\lim_{n\rightarrow \infty} \frac{\mathbb{E}\left[\Np{j}{k}\right]}{n} = \per\sum_{d^- = j}^\infty \sum_{d^+ = k}^\infty \diprob{d^-}{d^+}  {d^- \choose j}{d^+ \choose k}\per^{j+k}\left(1-\per\right)^{d^- - j + d^+ -k}  .
\end{align}
For $(j,k) = (0,0)$, we need to use equation \eqref{eq:def_exp_Dp_00} instead of equation \eqref{eq:def_exp_Dp}.
Since  $\Nr{d^-}{d^+}  \leq \N{d^-}{d^+}$ and  $\diarray$ is \proper\!\!,  for all $\epsilon > 0$ there exist $\kappa\left(\epsilon\right)$ and $N\left(\epsilon\right)$ such that for all $n > N$:
\begin{align*}
\sum_{\mathclap{\substack{(d^-,d^+) = (0,0)\\
			d^- \geq \kappa +1 \,\text{or}\, d^+ \geq \kappa +1}}}^{(d_{\max}, d_{\max})} \left(\N{d^-}{d^+} - \Nr{d^-}{d^+}\right)  \leq  \sum_{\mathclap{\substack{(d^-,d^+) = (0,0)\\
			d^- \geq \kappa +1 \,\text{or}\, d^+ \geq \kappa +1}}}^{(d_{\max}, d_{\max})} \N{d^-}{d^+} \leq \epsilon n.
\end{align*}
Thus the equivalent of \eqref{eq:bound_Dp_first} for $(j,k) = (0,0)$ becomes:
\begin{align*}
\sum_{ d^- =0}^{\kappa}\sum_{d^+ =0}^{ \kappa} \left[ \left(\N{d^-}{d^+} - \Nr{d^-}{d^+}\right )+P_{0,0}\left(d^-,d^+\right) \Nr{d^-}{d^+}\right] \leq \mathbb{E}\left[\Np{0}{0} \right] \leq \\  \sum_{ d^- =0}^{\kappa}\sum_{d^+ =0}^{ \kappa}\left[\left(\N{d^-}{d^+} - \Nr{d^-}{d^+}\right)  + P_{0,0}\left(d^-,d^+\right) \Nr{d^-}{d^+}\right] + 2\epsilon n,
\end{align*} 
and the analogous argument as for $(j,k) \neq (0,0)$ is applied to obtain:
\begin{align}
\label{eq:psite00}
\lim_{n\rightarrow \infty} \frac{\mathbb{E}\left[\Np{0}{0}\right]}{n} = \left(1-\per\right) + \per\sum_{\substack{d^- = j,\\d^+ = k}}^\infty\diprob{d^-}{d^+}  {d^- \choose j}{d^+ \choose k}\per^{j+k}\left(1-\per\right)^{d^- - j + d^+ -k}  = \psite{0}{0}.
\end{align}
Comparing equations \eqref{eq:psite} and \eqref{eq:psite00} with equation \eqref{eq:pbond} we obtain equation \eqref{eq:psite_pbond}.
\end{proof}
\end{Lemma}
Analysing equation \eqref{eq:psite_pbond}, we can see that  $\psite{j}{k}$ is correctly normalised and has moments:
\begin{align}
\label{eq:mum_site}
\mum^{\per, \text{site}}: =\mum[10]^{\per, \text{site}} =\mum[01]^{\per, \text{site}} = \per\mum^{\per, \text{bond}} = \per^2 \mum \quad \text{and} \quad \mum[11]^{\per, \text{site}} = \per \mum[11]^{\per, \text{bond}} = \per^3\mum[11].
\end{align}
Additionally, the generating function $U_{\per}^\text{site}(x,y)$ for $\psite{j}{k}$ is given by
\begin{equation}
\label{eq:gfs_Site}
\begin{aligned}
 U^\text{site}_{\per}(x,y) :=    1-\per+\per U(1-\per+\per x,1-\per+\per y),
 \end{aligned}
 \end{equation}
 and auxiliary functions $U^+_{\per},U^-_{\per}$ are the same as in equation \eqref{eq:gfs}.
It is left to determine $\hat{\per}^{\text{site}}$. If applicable, Theorem~\ref{thm:giant_scc} states that the percolation threshold is such $\per=\hat{\per}^{\text{site}}$ that 
$
\sum_{j,k=0}^\infty jk \psite{j}{k} = \sum_{j,k=0}^\infty j\psite{j}{k}.
$
Combing this  with equation \eqref{eq:mum_site} we find that the percolation thresholds for site and bond percolation are equal:
$$\hat{\per}^{\text{site}} = \frac{\mum}{\mum[11]} = \hat{\per}^{\text{bond}}.$$
This can be explained by noting that the expected degree distribution after site percolation is a rescaled version of the degree distribution after bond percolation, expect for $(0,0)$. Hence one expects the GSCC to appear under the same conditions. However, the GSCC after site percolation is expected to contain fewer vertices, because the probability to find an isolated vertex is larger.  
\subsubsection{Determining $\pcs$ and $\css$}
\label{subsubsc:theorem_site}
To finalise the proof of Theorems~\ref{thm:percolation_di_strong} and~\ref{thm:percolation_di_strong2} for site percolation, it remains to show that $\pcs = \hat{\per}^{\text{site}}$ and to determine $\css$.
This is done analogously to the proof for bond percolation  in Section~\ref{subsubsc:theorem_bond}. Because of the similarity between these proofs, we only explain the changes that are made in Section~\ref{subsubsc:theorem_bond} to convert it into the proof for site percolation.
 
 First, we need to replace $\pbond{j}{k}$  with $\psite{j}{k}$. 
Lemma~\ref{lm:site_cm} proves  the equivalent statement  for site percolation as Lemma~\ref{lm:percolation_is_CM} for bond percolation, therefore substituting 
  this lemma in Section~\ref{subsubsc:theorem_bond} will suffice. 
   However, equation  \eqref{eq:sum_prob_finite} requires a different proof. 
   Conditional on a certain realisation of $\left( \instubr, \outstubr\right)$ and the values $\si,\so \in I'_n$, $\ri,\ro \in I_n$, the value of $\Np{j}{k}$ is determined by the random choice of $\left(\instubm,\outstubm\right)$. By changing one element of $\left(\instubm,\outstubm\right)$ the value of $\Np{j}{k}$ changes by at most $2$. Thus Corollary~\ref{cor:a-b_sets} can be applied to obtain:
\begin{align*}
\p{ \left| \Np{j}{k} - \mathbb{E}\left[\Np{j}{k}\right] \right| > \sqrt{n}\ln^2 n \, \mid \,\si, \so, \ri, \ro, \left(\instubr, \outstubr\right)}\\\leq 2\exp \left( \frac{n \ln^2 n}{\left(m(1-\per)\per + n^{2/3} \ln^2 n\right)}\right) = e^ {-\Omega\left(\ln^2 n\right)}.
\end{align*}
Using Propositions~\ref{prop:I'}~and~\ref{lm:size_b-_b+}, it follows that 
\begin{align*}
\p{ \left| \Np{j}{k} - \mathbb{E}\left[\Np{j}{k}\right] \right| > \sqrt{n}\ln^2 n} = o\left(\frac{1}{n^3}\right),
\end{align*}
and, since $\kappa$ is bounded, this completes the proof of equation \eqref{eq:sum_prob_finite}.

The last
 change, is related to the fact that Theorem~\ref{thm:giant_scc} is now applied to a proper degree progression with $\diprobdissite$ as degree distribution instead of $\diprobdisbond$.   In Section~\ref{subsubsc:degree_dis_site} we found that  $\hat{\per}^{\text{site}} = \frac{\mum}{\mum[11]}$, and as in the case of bond percolation, we have  $U^+_{\per}(0),U^-_{\per}(0)>0$ for $0<\pi<1,$ which fulfils the last prerequisites for Theorem~\ref{thm:giant_scc}.
 Furthermore, in analogy to calculations in Section~\ref{subsubsc:theorem_bond}, 
we choose $x^*, \;y^*$ to be the solution of equation \eqref{eq:xystar} to find that:
 $$
 \css:=1-U_{\per}^\text{site}(x^*,1)-U_{\per}^\text{site}(1,y^*)+U_{\per}^\text{site}(x^*,y^*)=\per \csb.
$$
 which completes the proof of Theorems~\ref{thm:percolation_di_strong} and~\ref{thm:percolation_di_strong2} for site percolation. 
 \section*{Acknowledgments}
 The authors are grateful to Rik Versendaal for corrections and improvements to earlier drafts.
 

\end{document}